# Branching diffusions, superdiffusions and random media*

János Engländer

*Department of Statistics and Applied Probability
University of California, Santa Barbara, USA.
e-mail:* englander@pstat.ucsb.edu

**Abstract:** Spatial branching processes became increasingly popular in the past decades, not only because of their obvious connection to biology, but also because superprocesses are intimately related to nonlinear partial differential equations. Another hot topic in today's research in probability theory is 'random media', including the now classical problems on 'Brownian motion among obstacles' and the more recent 'random walks in random environment' and 'catalytic branching' models. These notes aim to give a gentle introduction into some topics in spatial branching processes and superprocesses in deterministic environments (sections 2-6) and in random media (sections 7-11).

Received September 2007.

## Contents

---

*This is an expository paper based on notes written for a graduate course at the University of Bath in July 2007. The course was generously funded by EPSRC (EP/E05448X/1). The author owes many thanks to S. Harris and A. Kyprianou as well as to the University of Bath for the invitation and for their hospitality during his visit.

# 1. Introduction

Our main purpose is to present some recent topics for branching processes and superprocesses, with special focus on random media. The choice of the topics of course reflects the author's interest in the field and is necessarily incomplete. There are several excellent expository books on superprocesses in general, for example by Dawson, Dynkin, Etheridge and Perkins and here we suffice with a subjective choice of some recent developments.

An outline of these notes is as follows. After defining the basic notions in Section 2, we turn to the connection between spatial branching processes and nonlinear partial differential equations (Section 3). In the following section (Section 4) we present some results on local extinction and local exponential growth, which is also a good excuse to discuss a version of the 'spine technology', a relatively recent (and, in the author's opinion, quite beautiful) tool appearing in the literature. The problem of the Law of large Numbers for spatial branching processes and superprocesses is discussed in Section 5. At the end of the 'deterministic environment part', we picked two further topics in superprocesses: the so-called compact support property and the 'polar decomposition' for superprocesses. This second one is somewhat related to the Law of large Numbers discussed in the previous section, because when dealing with LLN, one many times first decomposes a certain 'weighted' (or $h$-transformed) superprocess into its total mass ('radial' part) and a Fleming-Viot-type process ('angular' part).

In the second part (Sections 7-11) we will discuss spatial branching processes in random media, first giving a short review on some related classical problems (without branching), namely on Brownian motion among Poissonian obstacles and random walks in random environment (Section 7). In the next three sections we introduce branching processes and superprocesses among Poissonian obstacles: hard obstacles, 'mild' obstacles, and some generalizations (Sections 8, 9, and 10, respectively.) Finally, in the last section we review some basic catalytic branching problems.

References (including work cited as well as further bibliographic items) are given at the end of each section. These are not intended to be full bibliographies on the subjects, but they are sufficient to give the backgrounds for the particular theorems and problems discussed.

**Acknowledgement.** This paper is based on notes written for a graduate course at the University of Bath in July 2007. During the course several people, including A. Cox, S. Harris, A. Kyprianou, P. Mörters, J. Pardo Millan and A. Sakai, contributed with valuable remarks to the lectures. The author owes thanks to each one of them.

Several helpful suggestions by an anonymous referee are also gratefully acknowledged. Finally, the author is grateful to Prof. Aldous for inviting these lecture notes for publication in Probability Surveys.

## 2. Branching processes and superprocesses — basic notions

### *2.1. Branching diffusions*

Spatial branching processes have two main ingredients : the *spatial motion* and the *branching mechanism*. We start with describing the spatial motion.

Let $L$ be a second order elliptic (differential) operator on the Euclidean domain $D \subseteq \mathbb{R}^d$ of the form

$$L = \frac{1}{2} \sum_{i,j=1}^{d} a_{ij}(x) \frac{\mathrm{d}^2}{\mathrm{d}x_i \, \mathrm{d}x_j} + \sum_{i=1}^{d} b_i(x) \frac{\mathrm{d}}{\mathrm{d}x_i},$$

where $a_{ij}, b_i, \ i,j = 1...,d$, are in the class $C^{1,\eta}(D)$, $\eta \in (0,1]$ (i.e. their first order derivatives exist and locally Hölder-continuous), and the symmetric matrix $(a_{ij}(x))$ is positive definite on $D$.

Of course, $L$ can also be written in the slightly different 'divergence form':

$$L = \frac{1}{2} \nabla \cdot \widetilde{a} \nabla + \widetilde{b} \cdot \nabla.$$

The operator $L$ then corresponds to a *diffusion process* (or diffusion) on $D$ in the following sense. Take a sequence of increasing domains $D_n \uparrow D$ with $D_n \subset\subset D_{n+1}$[1] and let $\tau_{D_n}$ denote the first exit time from the (open) set $D_n$. Then there exists a unique family of probability measures $P_x$, $x \in D$ ($x$ is the starting point) describing the law of a Markov process $Y$ satisfying that

1. $P_x(Y_0 = x) = 1$,
2. $f(Y_{t \wedge \tau_{D_n}}) - \int_0^{t \wedge \tau_{D_n}} (Lf)(Y_s) \, ds$ is a martingale for all $f \in C^2(D)$ and all $n \geq 1$.

We say that the *generalized martingale problem* on $D$ has a unique solution and it is the law of the corresponding diffusion process or *L-diffusion* $Y$ on $D$.

Note that it is possible that the event $\lim_{n \to \infty} \tau_{D_n} < \infty$ ('explosion') has positive probability. In fact $\lim_{n \to \infty} \tau_{D_n} < \infty$ means that the process reaches $\Delta$, the 'cemetery state' in finite time, where $\Delta$ is identified with the Euclidean boundary of $D$ in the case of a bounded domain, or the Euclidean boundary of $D$ augmented by a point 'at infinity' in the case of an unbounded domain. In other words, the process actually lives on $\widehat{D} = D \cup \{\Delta\}$, the one-point compactification of $D$ and once it reaches $\Delta$ it stays there forever.

The connection to linear partial differential equations is well known. For a bounded continuous function $f$, consider the parabolic equation:

$$\dot{u} = Lu \text{ in } (0, \infty) \times D \tag{1}$$

with initial condition

$$u(0, \cdot) = f(\cdot) \text{ in } D.$$

---

[1] the notation $A \subset\subset B$ means that $A$ is bounded and its closure is still in $B$.

This Cauchy problem (generalized heat equation) is then solved by $u$, where $u(t,x) := E_x f(Y_t)$ and $Y$ is the diffusion corresponding to $L$ on $D$.

Furthermore, if $T_t f$ is defined by $T_t(f)(x) := E_x f(Y_t)$, then the Markov property of $Y$ yields that $\{T_t\}_{t \geq 0}$ is a semigroup, that is $T_{t+s} = T_t \circ T_s$ for $t, s \geq 0$. Now, (1) gives

$$\lim_{h \downarrow 0} \frac{T_{t+h} f(x) - T_t f(x)}{t} = L(T_t f)(x), \ t > 0$$

and formally we obtain ($t = 0$) that

$$\lim_{h \downarrow 0} \frac{T_h f - f}{t} = Lf,$$

which can in fact be verified for a certain class of $f$'s. Hence $L$ is often referred to as the *infinitesimal generator* of $Y$.

Diffusion processes behave in a very nice way from the point of view of their large time behavior, namely there are exactly two cases. Either

$$\forall x \in D, \emptyset \neq B \subset\subset D \text{ open}, \ P_x(Y_t \in B \text{ for arbitrarily large t's}) = 1$$

or

$$\forall x \in D, \emptyset \neq B \subset\subset D \text{ open}, \ P_x(Y_t \notin B \text{ for all } t > T(B, \omega)) = 1.$$

In the first case we say that $Y$ is *recurrent* and in the second we say that $Y$ is *transient*.

After describing the motion component ($L$-diffusion) we now turn to the branching component. Suppose that we start with a single individual and any individual has 0,1,2,... offspring with corresponding probabilities $p_0, p_1, p_2, ...$ and suppose that branching occurs at every time unit. If $h$ is the generating function of the offspring distribution,

$$h(z) := p_0 + p_1 z + p_2 z^2 + ...,$$

then clearly

$$h'(1) = p_1 + 2p_2 + 3p_3 + ... =: e,$$

and $e$ is the expected (or mean) offspring number.

As is well known, there are three cases for the Galton-Watson process defined above:

1. $h'(1) < 1$ and the system dies out in finite time (*subcritical* case)
2. $h'(1) = 1$ and the system dies out in finite time (*critical* case)
3. $h'(1) > 1$ and the system survives with positive probability (*supercritical* case)

For example we are talking about *strictly dyadic branching* when $p_2 = 1$ and in this case $h(z) = z^2$ and $h'(1) = 2$.

In the supercritical case $e_n$, the expected offspring number at time n, satisfies $e_n = e^n$ and when $p_0 > 0$, $d = P(\text{extinction})$ is the only root of $h(z) = z$ in $(0, 1)$. (The function $h$ can be shown to be concave upward.)

Many times, instead of unit time branching one considers random branching times with exponential distribution. In this case the above classification of subcritical, critical and supercritical branching is the same as before. Moreover, if the exponential rate is $\beta > 0$, then the probability of not branching up to time $t$ is obviously $e^{-\beta t}$ and a well known limit theorem due to Biggins states that if $Z_t$ is the population size at $t$ then in the supercritical case

$$\exists \lim_{t\to\infty} e^{-\beta t} Z_t = N > 0, \ a.s.,$$

on the survival set, as long as the offspring number $X$ satisfies that $X \log X$ has finite expectation. (There is a similar theorem for unit time branching too.)

After discussing motion and branching separately, let us try to put the two ingredients together! To this end, let $\beta \geq 0$ be smooth (more precisely, Holder continuous with some parameter $\eta \in (0,1]$) and not identically zero. The $(L,\beta)-$ *branching diffusion*[2] $X$ has the following (informal) definition.

**Definition 1** $((L,\beta)-$ branching diffusion)**.** Start with one initial particle located at some point $x \in D$, performing an $L$-diffusion on $D$. Her exponential clock depends on her position in the sense that the instantaneous rate of branching at $y \in D$ is $\beta(y)$. So for example the probability the particle has not branched by $t > 0$ is $\exp(-\int_0^t \beta(Y_s) \, ds)$. When the clock rings she is replaced by two offspring (strictly dyadic branching) [3]. The two offspring, starting from the position of the parent particle perform independent $L$-diffusions and are equipped with their independent exponential branching clocks, where again the instantaneous rate of branching at $y \in D$ is $\beta(y)$; etc. This way one obtains the stochastic process $Z$ (the value of $Z_t$, $t \geq 0$ is a particle configuration, or discrete measure on $D$); we will denote by $Z_t(B)$ the number of particles at $t$ in $B \subseteq D$.

Since $\beta$ is spatially dependent it is no longer clear what one means by subcritical, critical or supercritical case. Also, the criterion for extinction becomes a nontrivial question. One has to distinguish between extinction and *local extinction*.

**Definition 2.** $Z$ exhibits local extinction if for all $B \subset\subset D$, there exists an a.s. finite random time $\tau_B$ such that $Z_t(B) = 0$ for all $t > \tau_B$.

Intuitively, whether or not local extinction occurs should depend on the tradeoff between two properties:

1. how transient $Y$ (corresponding to $L$) is, and
2. how large $\beta$ is.

So, for example if $Y$ is 'very transient' comparing to the size of $\beta$, then the motion component 'wins' and $Z$ exhibits local extinction.

---

[2]Sometimes, when the domain is emphasized, one writes $(L, \beta, D)-$ branching diffusion.
[3]Obviously one could consider more generally a random number of offspring (possibly zero) according to any given branching law.

The 'strength of transience' can be measured by a single number $\lambda_c(L) \leq 0$ called the *generalized principal eigenvalue* of $L$ on $D$. We will give a different definition later but note here that for $B \subset\subset D$ open and $x \in B$,

$$\lim_{t \to \infty} \frac{1}{t} \log P_x(Y_s \in B, \text{ for all } s \leq t) = \lambda_B < 0,$$

where $\lambda_B$ does not depend on $x \in B$, and $\lambda_{B_n} \uparrow \lambda_c(L)$ for a sequence of domains $B_n \uparrow D$. So $\lambda_c$ is loosely speaking the 'exponential rate of escape from compacts'.

The main result concerning local extinction (we will discuss its more general form later) is that when $\beta$ is a positive constant, local extinction occurs if and only if $\lambda_c(L) \leq -\beta$.

We close this section with two definitions of superprocesses — as scaling limits of branching diffusions[4] and also as measure valued Markov processes via their Laplace functionals. We start with the latter one.

### 2.2. Superprocess via its Laplace functional

Let $D \subseteq \mathbb{R}^d$ be a domain and let $\mathcal{B}(D)$ denote the Borel sets of $D$. We write $\mathcal{M}_f(D)$ and $\mathcal{M}_c(D)$ for the class of finite measures resp. the class of finite measures with compact support on $\mathcal{B}(D)$. For $\mu \in \mathcal{M}_f(D)$, denote $\|\mu\| := \mu(D)$ and let $C_b^+(D)$ and $C_c^+(D)$ be the class of non-negative bounded continuous resp. non-negative bounded continuous functions $D \to \mathbb{R}$ having compact support. Write $C^{k,\eta}(D)$ for the usual Hölder spaces of index $\eta \in (0,1]$ including derivatives of order $k$, and set $C^\eta(D) := C^{0,\eta}(D)$.

Let $L$ be given (in divergence form) by

$$L := \frac{1}{2} \nabla \cdot a \nabla + b \cdot \nabla, \qquad (2)$$

where $a_{i,j}, b_i \in C^{1,\eta}(D)$, $i,j = 1,...,d$, for some $\eta \in (0,1]$, and the matrix $a(x) := (a_{i,j}(x))$ is symmetric, and positive definite for all $x \in D$. In addition, let $\alpha, \beta \in C^\eta(D)$, and assume that $\alpha$ is positive, and $\beta$ is bounded from above.

**Definition 3** (($L, \beta, \alpha; D$)-superdiffusion)**.** With $D, L, \beta$ and $\alpha$ as above, let $(X, \mathbf{P}_\mu, \mu \in \mathcal{M}_f(D))$ denote the $(L, \beta, \alpha; D)$-superdiffusion, where $\mu$ denotes the starting measure $X_0$. That is, $X$ is the unique $\mathcal{M}_f(D)$-valued continuous (time-homogeneous) Markov process which satisfies, for any $g \in C_b^+(D)$,

$$\mathbf{E}_\mu \exp \langle X_t, -g \rangle = \exp \langle \mu, -u(\cdot, t) \rangle, \qquad (3)$$

where $u$ is the minimal nonnegative solution to

$$\left. \begin{array}{c} u_t = Lu + \beta u - \alpha u^2 \quad \text{on } D \times (0, \infty), \\ \lim_{t \downarrow 0} u(\cdot, t) = g(\cdot). \end{array} \right\} \qquad (4)$$

As usual, $\langle \nu, f \rangle$ denotes the integral $\int_D \nu(\mathrm{d}x) f(x)$.

---

[4]As already mentioned, branching diffusions are also measure valued processes: $Z_t(\cdot) = \sum_{i=1}^{N_t} \delta_{Z_t^i}(\cdot)$, where $N_t$ is the population size at time $t$.

One usually refers to $\beta$ as *mass creation* and $\alpha$ as the *intensity parameter* (or variance). Finally, the Markov property is in fact equivalent to the property that the shift defined by (4) defines a semigroup, which in turn follows from the minimality of the solution.

### 2.3. The particle picture for the superprocess

Previously we defined the $(L, \beta, \alpha; D)$-superprocess $X$ analytically. In fact, $X$ also arises as the short life time and high density diffusion limit of a *branching particle system*, which can be described as follows: in the $n^{\text{th}}$ approximation step each particle has mass $1/n$ and lives a random time which is exponential with mean $1/n$. While a particle is alive, its motion is described by a diffusion process corresponding to the operator $L$ (on $\widehat{D} = D \cup \{\Delta\}$). At the end of its life, the particle located at $x \in D$ dies and is replaced by a random number of particles situated at the parent particle's final position. The law of the number of descendants is spatially varying such that the mean number of descendants is $1 + \frac{\beta(x)}{n}$, while the variance is assumed to be $2\alpha(x)$. All these mechanisms are independent of each other.

More precisely, for each positive integer n, consider $N_n$ particles, each of mass $\frac{1}{n}$, starting at points $x_i^{(n)} \in D, i = 1, 2, \ldots, N_n$, and performing independent branching diffusion according to the motion process $Y$, with branching rate $cn, c > 0$, and branching distribution $\{p_k^{(n)}(x)\}_{k=0}^{\infty}$, where

$$e_n(x) \equiv \sum_{k=0}^{\infty} k p_k^{(n)}(x) = 1 + \frac{\gamma(x)}{n},$$

and

$$v_n^2(x) \equiv \sum_{k=0}^{\infty} (k-1)^2 p_k^{(n)}(x) = m(x) + o(1)$$

as $n \to \infty$, uniformly in $x$; $m, \gamma \in C^{\eta}(D), \eta \in (0, 1]$ and $m(x) > 0$. Let

$$\mu_n = \frac{1}{n} \sum_{i=1}^{N_n} \delta_{x_i^{(n)}}.$$

Let $N_n(t)$ denote the number of particles alive at time t and denote their positions by $\{X_i^n(t)\}_{i=1}^{N_n(t)}$. Denote by $\mathcal{M}_f(D)$ ($\mathcal{M}_f(\widehat{D})$) the space of finite measures on $D$ ($\widehat{D}$). Define an $\mathcal{M}_f(\widehat{D})$-valued process $X_n$ by

$$X_n(t) = \frac{1}{n} \sum_{i=1}^{N_n(t)} \delta_{X_i^n(t)}.$$

Denote by $P_{\mu_n}^{(n)}$ the probability measure on $D([0, \infty), \mathcal{M}_f(\widehat{D}))$ induced by $X_n$. Assume that $m$ and $\gamma$ are bounded from above. Let $w - \lim_{n \to \infty} \mu_n = \mu \in$ $\mathcal{M}_f(D)$. (The notation $w - lim$ denotes the limit in the weak topology.) Then there exists a $P_\mu^* \in C([0,\infty), \mathcal{M}_f(\widehat{D}))$ such that $P_\mu^* = w - \lim_{n\to\infty} P_{\mu_n}^{(n)}$. Define $P_\mu \in C([0,\infty), \mathcal{M}_f(D))$ by $P_\mu(\cdot) = P_\mu^*(\cdot \cap D)$ and let $X$ be the process corresponding to $P_\mu$. Then $X$ is an $(L, \beta, \alpha; D)$-superprocess, where $L$ corresponds to $Y$ on $D$, $\beta(x) = c\gamma(x)$ and $\alpha(x) = \frac{1}{2}cm(x)$.

Finally, in the particular case when $\beta \equiv 0$, an alternative approximation yields the same superprocess: in the $n^{\text{th}}$ approximation step one considers critical branching diffusions (the motion component corresponds to $L$ and the branching is strictly binary, i.e. either zero or two offspring), but the branching rate is now $2n\alpha(x)$. So $\alpha(\cdot)$ in this case can also be thought of as the branching 'clock'.

## 3. Connection between spatial branching processes and nonlinear partial differential equations; the qualitative behavior of superprocesses

### *3.1. Log-Laplace equations*

Let us recall the log-Laplace equation for superprocesses: for any $g \in C_b^+(D)$,

$$\mathbf{E}_\mu \exp \langle X_t, -g \rangle = \exp \langle \mu, -u(\cdot, t) \rangle,$$

where $u$ is the minimal nonnegative solution to (4). This fact already shows that there exists a connection between measure valued processes and partial differential equations.

**Remark 1.** A similar connection with the $(L, \beta)$-branching diffusion $Z$ (with probability $P$) is as follows. For any $g \in C_b^+(D)$, $u(x,t) = E_x \exp \langle Z_t, -g \rangle$ solves

$$\left. \begin{array}{c} u_t = Lu + \beta(u^2 - u) \quad \text{on } D \times (0, \infty), \\ \lim_{t \downarrow 0} u(\cdot, t) = e^{-g(\cdot)}, \\ 0 \leq u \leq 1. \end{array} \right\} \qquad (5)$$

Note that, unlike in the superprocess case, here the $u^2$ term simply comes from the fact that the branching is strictly dyadic and so the generating function is just $h(z) = z^2$. Note also that in (5) one can always 'switch' from $u^2 - u$ to $v - v^2$ as $v := 1 - u$ solves

$$\left. \begin{array}{c} v_t = Lv + \beta(v - v^2) \quad \text{on } D \times (0, \infty), \\ \lim_{t \downarrow 0} v(\cdot, t) = 1 - e^{-g(\cdot)}, \\ 0 \leq v \leq 1. \end{array} \right\}$$

$\diamond$

### *3.2. Exit measure and Brownian snake*

Dynkin's *exit measure*[5] yields yet another connection with semilinear partial differential equations. Informally, the exit measure $X^D$ from a domain $D \subset \mathbb{R}^d$ is obtained by 'freezing' mass of the superprocess when it first hits the boundary $\partial D$. One way of defining the exit measure is through the approximating particle picture: the exit measure can easily be defined for the $n$th level discrete particle system, and then one has to verify that as $n \to \infty$, those discrete exit measures converge to a limit. (See [5] for more on exit measures.)

Although in this subsection we will talk about Brownian motion, in fact the method can be extended from Brownian motion (and equations involving $\Delta$) to general diffusions (and equations with elliptic operators $L$).

---

[5]We will revisit a version of it in Subsection 4.1.

Dynkin showed that if $D$ is smooth and $g$ is a bounded function on $\partial D$, then $u(x) = -\log \mathbb{E}_{\delta x} \exp -\langle X^D, g \rangle, x \in D$ solves the problem

$$\left. \begin{aligned} \Delta u &= u^2 \quad on\ D, \\ u|_{\partial D} &= g. \end{aligned} \right\}$$

In order to investigate solutions with *boundary blow-up*, choose $g_n \equiv n$ and then let $n \uparrow \infty$, quickly leading to the following representation.[6]

**Corollary 1.** *Let $D \subset \mathbb{R}^d$ be a bounded regular [7] domain. The function $u(x) = -\log \mathbb{P}_{\delta x}(X^D = 0), x \in D$ is the minimal nonnegative solution to the problem*

$$\left. \begin{aligned} \Delta u &= u^2 \quad on\ D, \\ u|_{\partial D} &= +\infty. \end{aligned} \right\}$$

Before turning to the next part, let us quickly mention that Dynkin actually defined the exit measure not only from spatial but also from *space-time* domains (c.f. Subsection 4.1), leading to probabilistic solutions to nonlinear *parabolic* boundary value problems. (We do not discuss this further here, the interested reader should consult [5].)

Le Gall introduced a path-valued process called the *Brownian snake* and related it to the exit measure, thus providing yet another probabilistic representation for nonlinear PDE's (in terms of the snake).

The Brownian snake combines the genealogical structure of random real trees with spatial motions governed by a general Markov process (usually taken as Brownian motion). Informally, each Brownian snake path corresponds to the spatial positions along the ancestral line of a vertex in the branching tree. The precise definition of the Brownian snake is motivated by the 'coding' of real trees. In view of applications to PDE, Le Gall characterized the exit measure as a measure, which, in a sense, is uniformly spread over the set of exit points of the Brownian snake paths from $D$. (See [9] for the precise definition of the snake and the characterization of the exit measure.)

When the underlying spatial motion is $d$-dimensional Brownian motion, the connection between the exit measure and the semilinear PDE $\Delta u = u^2$ is as follows.

Recalling that the snake is a path valued process and that the exit measure is expressed in terms of the snake, Le Gall defines a $\sigma$-finite measure on $C(\mathbb{R}_+, \mathcal{W})$ denoted by $\mathbb{N}_x$[8], where $\mathcal{W}$ is the set of all finite paths (see [9] for the precise definition of $\mathbb{N}_x$), and proves the following representation (notice the extra factor 4 in the equation).

---

[6] In [12, 8] there is a different approach to solutions with boundary blow-up without using the exit measure.

[7] The domain $D$ is called regular if every $x \in \partial D$ satisfies that $\inf\{t > 0 \mid \xi_t \notin D\} = 0$ $P_x$-a.s., where $(\xi, P)$ is Brownian motion

[8] It is dubbed the 'excursion measure away from $x$'.

**Proposition 1.** *Let $D$ be a bounded regular domain. Let $g$ be a continuous nonnegative function on $\partial D$. Then $u(x) = \mathbb{N}_x(1 - \exp-\langle X^D, g\rangle), x \in D$ solves the problem*

$$\left.\begin{aligned}\Delta u &= 4u^2 \quad \text{on } D,\\ u|_{\partial D} &= g.\end{aligned}\right\}$$

Similarly as before, this leads to the following representation of the solution with boundary blow-up.

**Corollary 2.** *Let $D$ be a bounded regular domain. Then the function $u(x) = \mathbb{N}_x(X^D \neq 0), x \in D$ is the minimal nonnegative solution of the problem*

$$\left.\begin{aligned}\Delta u &= 4u^2 \quad \text{on } D,\\ u|_{\partial D} &= +\infty.\end{aligned}\right\}$$

Due to groundbreaking work of Dynkin, Le Gall and Mselati (see [5, 9, 11]), we now have a good understanding of the classification and probabilistic representation of solutions of $\Delta u = u^2$ in a domain.

### 3.3. Support properties

It turns out that (4) can be used very efficiently in investigating the following properties of the superprocess.

**Definition 4.** A path $X.$ *survives* if $X_t \neq 0$, $\forall t \geq 0$ and becomes *extinct* otherwise.

We will use the shorthand $S := \{X_t \neq 0, \ \forall t \geq 0\}$ for the event of survival. The basic question is to calculate $\mathbf{P}_\mu(S)$, and in particular to decide whether this probability is zero ('extinction').

**Definition 5.** The superprocess exhibits *local extinction* if for every set $B \subset\subset D$ there exists a $\mathbf{P}_\mu$-a.s. finite random time $\tau_B$ such that $X_t(B) = 0$, $\forall t > \tau_B$.

Unlike for discrete particle systems, however, it is not clear whether '$B$ becomes eventually empty' is the same as $\lim_{t\to\infty} X_t(B) = 0$. R. Pinsky [12] introduced the following notions.

**Definition 6.** Assume that $\mathbf{P}_\mu(S) > 0$.

- The support of $X$ is *recurrent* if

$$\mathbf{P}_\mu(X_t(B) > 0, \text{ for some } t \mid S) = 1, \ \forall B \subset D \text{ open}.$$

- When $d \geq 2$, the support of $X$ is *transient* if

$$\mathbf{P}_\mu(X_t(B) > 0, \text{ for some } t \mid S) < 1, \ \forall B \subset\subset D \text{ s.t.}$$

$D \setminus B$ is connected and $\mathrm{supp}(\mu) \cap \overline{B} = \emptyset$. [$\overline{B}$ is the closure of $B$.]

- When $d = 1$, the support of $X$ is *transient* if

$$\mathbf{P}_\mu(X_t(B) > 0, \text{ for some } t \mid S) < 1,$$

either for all $B \subset\subset D$ s.t. $\inf \text{supp}(\mu) > \sup B$ or for all $B \subset\subset D$ s.t. $\sup \text{supp}(\mu) < \inf B$.

It turns out that in the above definitions there is no dependence on the starting measure, i.e. each property (extinction, survival, local extinction, recurrence and transience) are either true for all $\emptyset \neq \mu \in \mathcal{M}_c$ or not true for any $\emptyset \neq \mu \in \mathcal{M}_c$.

### 3.4. Weighted superprocess and nonlinear h-transform

A very important tool in the investigations of these *support properties* is a transformation that leaves the support invariant but changes the mass. If $X$ is a superprocess and $h > 0$ is smooth, then we will define below a *weighted superprocess* $X^h$ which has the same support but has different motion and branching mechanism. Let

$$X^h := hX, \quad \text{i.e } \frac{dX_t^h}{dX_t} = h, \ t \geq 0,$$

and notice that the log-Laplace equation with $\phi = h\psi$, $\mu = \nu/h$ yields

$$\mathbf{E}_\nu e^{\langle -\psi, X_t^h \rangle} = \mathbf{E}_\mu e^{\langle -\phi, X_t \rangle} = e^{\langle -u(\cdot, t), \mu \rangle},$$

where $u$ is the minimal nonnegative solution to (4) with $\phi$ in place of $g$.

It is sufficient to restrict the setting on compactly supported measures, $\nu \in \mathcal{M}_c \Leftrightarrow \mu \in \mathcal{M}_c$ and test functions, $\phi \in C_c^+(D) \Leftrightarrow \psi \in C_c^+(D)$.

Let $\mathcal{A}(u) := Lu + \beta u - \alpha u^2$. Now, if $v := u/h$, then

$$\mathbf{E}_\nu e^{\langle -\psi, X_t^h \rangle} = e^{\langle -v(\cdot, t), \nu \rangle},$$

where

$$v_t = \frac{1}{h}\mathcal{A}(hv)$$

and $v(\cdot, 0) = \psi(\cdot)$. What this means is that $X^h$ is another superprocess! If we now introduce the notation $\mathcal{A}^h(\cdot) := \frac{1}{h}\mathcal{A}(h \cdot)$ (analogously to the Doob's $h$-transform for linear operators), then we obtain the log-Laplace equation for the weighted superprocess $X^h$ in the following form: for any $\psi \in C_c^+(D)$,

$$\mathbf{E}_\nu \exp \langle X_t, -g \rangle = \exp \langle \nu, -v(\cdot, t) \rangle, \tag{6}$$

where $v$ is the minimal nonnegative solution to

$$\left.\begin{aligned} v_t &= \mathcal{A}^h(v) \quad \text{on } D \times (0, \infty), \\ \lim_{t \downarrow 0} u(\cdot, t) &= \psi(\cdot). \end{aligned}\right\} \tag{7}$$

That is, if $X$ corresponds to the operator $\mathcal{A}$ then $X^h$ corresponds to the operator $\mathcal{A}^h$. Although $X^h$ is just a weighted version of $X$, it corresponds to a different motion and a different branching mechanism, since, as a straightforward computation reveals,

$$\mathcal{A}^h(v) = \left(L + a\frac{\nabla h}{h} \cdot \nabla + \frac{(L+\beta)h}{h}\right)v - \alpha h v^2.$$

So the new motion corresponds to $L + a\frac{\nabla h}{h} \cdot \nabla$ and the new branching mechanism corresponds to the term $\frac{(L+\beta)h}{h}v - \alpha h v^2$.

Just like the transformation $X \to X^h$ preserves the support of $X$, the corresponding transformation $\mathcal{A} \to \mathcal{A}^h$ (nonlinear $h$-transform) preserves some important properties of the operator. For example, the *generalized principal eigenvalue* of the linear part of the operator is invariant under the transformation.

### 3.5. Some results

In the previous section we discussed the generalized principal eigenvalue of a diffusion operator and its connection to local extinction for branching diffusions (with constant branching rate). We now discuss it in more generality.

Define the generalized principal eigenvalue as

$$\lambda_c = \lambda_c(L + \beta) := \inf\{\lambda \mid \exists u > 0 \ s.t. \ (L + \beta - \lambda)u = 0,\}$$

and note that for constant $\beta$, one has $\lambda_c(L+\beta) = \lambda_c(L) + \beta$. Now a more general result concerning local extinction is as follows.

**Theorem 1.** *The $(L, \beta, \alpha, D)$-superprocess (or the $(L, \beta, D)$-branching diffusion) exhibits local extinction if and only if $\lambda_c \leq 0$.*

(The superprocess version of this theorem is due to R. Pinsky.) For superprocesses, since extinction implies local extinction, one obtains that a sufficient condition for survival w.p.p. is that $\lambda_c(L) > 0$.

It is beyond the scope of these notes to discuss all the support properties, but as an illustration, we show how one can handle the extinction probability by using the log-Laplace equation.

First, notice that for $\mu \in \mathcal{M}_c$,

$$\mathbf{P}_\mu(X_t = 0) = \mathbf{E}_\mu \lim_{n\to\infty} e^{\langle -n, X_t\rangle} = \lim_{n\to\infty} \mathbf{E}_\mu e^{\langle -n, X_t\rangle} = \lim_{n\to\infty} e^{\langle -u^{(n)}(\cdot,t),\mu\rangle},$$

where $u^{(n)}$ is the minimal nonnegative solution to (4) with $g$ replaced by $n$. It is clear probabilistically that $u^{(n)}$ is monotone increasing in $n$ and it can be shown that $u := \lim_{n\to\infty} u^{(n)}$ is finite and it is the minimal nonnegative solution to (4) with $g$ replaced by $\infty$.

We have obtained
$$\mathbf{P}_\mu(X_t = 0) = e^{\langle -u(\cdot,t),\mu\rangle},$$

and now letting $t \uparrow \infty$ and defining $w(\cdot) = \lim_{t \to \infty} u(\cdot, t)$ (again, monotonicity is clear probabilistically), we obtain the probability of extinction:

$$\mathbf{P}_\mu(S^c) = e^{\langle -w, \mu \rangle}.$$

As expected, it can be shown that $w$ solves the steady-state equation $Lw + \beta w - \alpha w^2 = 0$. The question is now: which (nonnegative) solution of the steady-state equation is $w$? This is a subtle issue as $w$ is sometimes the *maximal* solution and sometimes it isn't.

More precisely, it can be shown that the maximal solution $w_{max}$ always exists and

$$\mathbf{P}_\mu(\mathcal{R} \subset\subset D) = e^{\langle -w_{max}, \mu \rangle},$$

where $\mathcal{R} := \overline{\bigcup_{s \geq 0} \text{supp}(X_s)}$ is the *range* of the process.

Since $w \leq w_{max}$, $\mathbf{P}_\mu(\mathcal{R} \subset\subset D) \leq \mathbf{P}_\mu(S^c)$.

Suppose now that the superprocess possesses the *compact support property*, that is, if $\mu$ is compactly supported then so is $\overline{\bigcup_{0 \leq s \leq t} \text{supp}(X_s)}$ for all $t > 0$ (this property too is independent of the starting measure). Then, obviously $\mathbf{P}_\mu(\mathcal{R} \subset\subset D) \geq \mathbf{P}_\mu(S^c)$, that is $\mathbf{P}_\mu(\mathcal{R} \subset\subset D) = \mathbf{P}_\mu(S^c)$ and $w = w_{max}$.

This observation leads to a nice way of proving explosion (no compact support property) for the superprocess. Namely, all one has to show is that $w$ is *not* the maximal nonnegative solution to the steady-state equation. For example, when we know that the process becomes extinct (and thus $w = 0$), the existence of a positive solution yields explosion.

Based on these ideas one can come up with examples where the underlying process is conservative (i.e. no explosion occurs), but nevertheless the superprocess explodes. This can be shown [7] to be true even in the case of a time changed one dimensional Brownian motion with $\beta = 0$, $\alpha = 1$ (i.e. in the case of a $(\rho(x)\Delta u, 0, 1; \mathbb{R})$-superprocess), if the time change $\rho > 0$ is sufficiently large. (Note that the time change does not effect recurrence and so the time changed Brownian motion is still conservative.)

Of course one may consider the $(\rho(x)\Delta u, 0, 1; \mathbb{R}^d)$-superprocess in any dimension $d \geq 1$, and wonder, how large $\rho$ has to be in order to loose the compact support property. Somewhat surprisingly, it turns out [7] that while the threshold for $\rho$ is quadratic for $d \geq 2$, it is cubic for $d = 1$.

## 4. Local extinction versus local exponential growth; the 'spine'

In this section we present a 'spine' construction and its application in the study of local extinction versus local exponential growth of spatial branching processes. Although it can be applied to superprocesses too, we will only discuss, the more intuitive[9] branching diffusion case.

Let $X$ denote the $(L,\beta)$-branching diffusion in the section. The finite (discrete) starting measure will be denoted by $\mu$.

### 4.1. The spine decomposition

**A 'natural martingale'.** Let $B \subset\subset D$ be a nonempty domain. Suppose that $\text{supp}(\mu) \subset B$. Let $\mathcal{F}_t$ denote the natural filtration of $X$ up to time $t \geq 0$ and let $X^{t,B}$ denote *Dynkin's exit measure* from $B \times [0,t)$. Loosely speaking, $X^{t,B}$ is obtained by 'freezing' particles on $\partial B$ at the time they hit it and so $X^{t,B}$ is a discrete measure on $(B \times \{t\}) \cup (\partial B \times (0,t))$.

Let $\lambda$ denote the generalized principal eigenvalue[10] of $L + \beta$ on $B$. It is a standard fact that there exists a unique nonnegative function $\phi$ which solves $(L + \beta - \lambda)\phi = 0$ on $B$ with Dirichlet boundary condition. Define for each fixed $t \geq 0$, $\phi^t : \overline{B} \times [0,t] \to [0,\infty)$ such that $\phi^t(\cdot, u) = \phi(\cdot)$ for each $u \in [0,t]$. (Here $\overline{B}$ denotes the closure of $B$). We claim that

$$\left\{ M_t^\phi := e^{-\lambda t} \frac{\langle \phi^t, X^{t,B} \rangle}{\langle \phi, \mu \rangle} : t \geq 0 \right\}$$

is a (mean one) $P_\mu$-martingale. To see this note that on account of the so-called Dynkin's Markov property of exit measures and the Dirichlet boundary condition for $\phi$,

$$E_\mu(M_{t+s}^\phi \mid \mathcal{F}_t) =$$

$$E_\mu \left( e^{-\lambda(t+s)} \frac{\langle \phi^{t+s}, X^{t+s,B} \rangle}{\langle \phi, \mu \rangle} \,\bigg|\, \mathcal{F}_t \right) = e^{-\lambda t} E_{X_t} \left( e^{-\lambda s} \frac{\langle \phi^s, X^{s,B} \rangle}{\langle \phi, \mu \rangle} \right)$$

for $s, t \geq 0$. We have to show that the righthand side equals $M_t^\phi$ a.s. On account of the branching property it is enough to show that $E_{\delta_x}(e^{-\lambda t} \langle \phi^t, X^{t,B} \rangle) = \phi(x)$ for all $t \geq 0$ and $x \in B$. To show this latter property, note that from the stochastic representation for backward equations: $v(x,0) := E_{\delta_x}(e^{-\lambda t} \langle \phi^t, X^{t,B} \rangle)$ where $v$ solves $-\dot{v} = (L + \beta - \lambda)v$ on $B \times (0,t)$ with $v = \phi^t$ on $\partial^{t,B}$ and $\partial^{t,B}$ is the union of $\partial B \times (0,t]$ and $B \times \{t\}$. The boundary condition follows from the regularity properties of the underlying diffusion. Therefore, by parabolic uniqueness, we obtain that $v(\cdot, 0) = \phi(\cdot)$, $t \geq 0$.

Our goal is to define a *change of measure* with this 'natural' martingale. (For the reader unfamiliar with this notion: if the new measure $\mathbb{Q}$ for a continuous

---

[9] The author's wife is not a mathematician. After listening to phone conversations with A. Kyprianou she asked 'who this Martin Gale was' and 'what was wrong with his spine'.

[10] which is just the usual principal Dirichlet eigenvalue.

time stochastic process with filtration $\{\mathcal{G}_t\}$ is defined by

$$\left.\frac{d\mathbb{Q}_x}{d\mathbb{P}_x}\right|_{\mathcal{G}_t} = M_t, \tag{8}$$

then it is easy to show that $M$ must be a mean one $\mathbb{P}_x$-martingale. Conversely, every mean one martingale defines a proper law for the stochastic process with the definition under (8).) Before we do that, let us first define two changes of measures: the first is for the motion component and the second is for a Poisson point process.

As before $B$ will always denote a nonempty open set compactly embedded in $D$ with a smooth boundary.

**Girsanov change of measure.** Let $\lambda = \lambda_c(L+\beta, B)$. Just like before, $\phi$ denotes the eigenfunction satisfying

$$(L+\beta-\lambda)\phi = 0 \text{ in } B \text{ with } \phi = 0 \text{ on } \partial B.$$

Let $\tau^B = \inf\{t \geq 0 : Y_t \notin B\}$ and assume that the diffusion $(Y, \mathbb{P}_x)$ is adapted to some filtration $\{\mathcal{G}_t : t \geq 0\}$. Then under the change of measure

$$\left.\frac{d\mathbb{P}_x^\phi}{d\mathbb{P}_x}\right|_{\mathcal{G}_t} = \frac{\phi(Y_{t\wedge\tau^B})}{\phi(x)} \exp\left\{-\int_0^{t\wedge\tau^B} \lambda - \beta(Y_s)\, ds\right\}$$

the process $(Y, \mathbb{P}_x^\phi)$ corresponds to the $h$-transformed ($h = \phi$) generator $(L+\beta-\lambda)^\phi = L + a\phi^{-1}\nabla\phi\cdot\nabla$.

It can be shown that the process $(Y, \mathbb{P}_x^\phi)$ is *ergodic* on $B$ (i.e. it is positive recurrent).

**Change of measure for Poisson point processes.** Given a nonnegative bounded continuous function $g(t)$, $t \geq 0$, the Poisson point process $\eta$ is defined as follows. Let $\{n = n_t; t \geq 0\}$ be a Poisson process with instantaneous rate $g(t)$, and let $\eta = \{\{\sigma_i : i = 1, ..., n_t\} : t \geq 0\}$ be the corresponding Poisson point process on $[0, \infty)$ ($\sigma_i$ is the $i$th 'arrival'). Let $\mathbb{L}^g$ denote the law of $\eta$.

If $\eta$ is adapted to $\{\mathcal{G}_t : t \geq 0\}$, then

$$\left.\frac{d\mathbb{L}^{2g}}{d\mathbb{L}^g}\right|_{\mathcal{G}_t} = 2^{n_t} \exp\left\{-\int_0^t g(s)\, ds\right\}.$$

Combining the change of measure for the motion component and for the Poisson point process leads to a change of measure for the spatial branching process $X$. Our 'spine' result here is as follows.

**Theorem 2.** *Suppose that $\mu$ is a finite measure with* $\operatorname{supp}\mu \subset B$. *For branching particle process we can thus take $\mu = \sum_i \delta_{x_i}$ where $\{x_i\}$ is a finite set of (not necessarily distinct) points in $B$. Define $\widetilde{P}_\mu$ by the martingale change of measure*

$$\left.\frac{d\widetilde{P}_\mu}{dP_\mu}\right|_{\mathcal{F}_t} = M_t^\phi.$$

*Define*

$$\mathbb{P}^\phi_{\phi\mu} = \int_B \mu\,(dx)\,\frac{\phi(x)}{\langle\phi,\mu\rangle}\,\mathbb{P}^\phi_x,$$

*that is, we randomize the starting point of* $(Y, \mathbb{P}^\phi_\cdot)$ *according to the probability distribution* $\phi\mu/\langle\phi,\mu\rangle$. *Note in particular that when* $\mu = \delta_x$, $\mathbb{P}^\phi_{\phi\mu} = \mathbb{P}^\phi_x$.

Suppose that $g \in C_b^+(D)$ and $u_g$ is the minimal positive solution to $\dot{u} = Lu + \beta u^2 - \beta u$ on $D \times (0,\infty)$ with $\lim_{t\downarrow 0} u(\cdot,t) = g(\cdot)$. Then

$$\widetilde{E}_\mu\left(e^{-\langle g, X_t\rangle}\right) =$$
$$\sum_i \frac{\phi(x_i)}{\langle\phi,\mu\rangle}\left\{\mathbb{E}^\phi_{x_i}\mathbb{L}^{2\beta(Y)}(e^{-g(Y_t)}\prod_{k=1}^{n_t} u_g(Y_{\sigma_k}, t-\sigma_k))\prod_{j\neq i} u_g(x_j, t)\right\} \quad (9)$$

The decomposition suggest that $(X, \widetilde{P}_\mu)$ has the same law as a process constructed in the following way. From the configuration $\mu = \sum_i \delta_{x_i}$ pick a point $x' \in \{x_i\}$ with probability $\phi(x')/\langle\phi,\mu\rangle$. From the remaining points, independent $(L,\beta;D)$-branching processes are initiated. From the chosen point, a $(Y, \mathbb{P}^\phi_{x'})$-diffusion is initiated along which $(L,\beta;D)$-branching processes immigrate at space-time points $\{(Y_{\sigma_i}, \sigma_i) : i \geq 1\}$ where $n = \{\{\sigma_i : i = 1,...,n_t\} : t \geq 0\}$ is a Poisson process with law $\mathbb{L}^{2\beta(Y)}$.

This decomposition relates to a spectrum of similar results that exist in the literature for both superprocesses and branching processes by S. Evans, A. Etheridge, R. Williams, S. Roelly, A. Rouault and others [11]. Theorem 2 offers decompositions with the particular feature that the spine is represented by a diffusion conditioned to stay in the compactly embedded domain $B$.

Now, obviously the new measure is absolutely continuous with respect to the old one *up to time t*. However this does not mean that this is true 'up to time infinity'. Nevertheless, we have the following result.

**Lemma 1.** *Suppose that* $\operatorname{supp}\mu \subset B$ *and* $\lambda = \lambda_c(L+\beta, B) > 0$. *Then* $M_t^\phi$ *converges to its almost sure limit* $M_\infty^\phi$ *in* $L^1(P_\mu)$, *and furthermore* $\widetilde{P}_\mu \ll P_\mu$.

### *4.2. Application to local extinction and exponential growth*

After all these preparations, let us see how one can attack the problem of local extinction and (local) exponential growth for the branching diffusion.

**Theorem 3** (Local extinction vs. local exponential growth). *Let* $0 \neq \mu$ *be a measure with* $\operatorname{supp}\mu \subset\subset D$.

(i) *Under* $P_\mu$ *the process* $Z$ *exhibits local extinction if and only if there exists a function* $h > 0$ *satisfying* $(L+\beta)h = 0$ *on* $D$, *that is, if and only if* $\lambda_c \leq 0$. *In particular, the local extinction property does not depend on the starting measure* $\mu$.

---

[11]The spine decomposition presented in this notes was proved by A. Kyprianou and the author and to a large extent it was inspired by a result of R. Lyons, R. Pemantle and Y. Peres.

(ii) *When $\lambda_c > 0$, for any $\lambda < \lambda_c$ and any open $\emptyset \neq B \subset\subset D$,*

$$P_\mu \left( \limsup_{t \uparrow \infty} e^{-\lambda t} Z_t(B) = \infty \right) > 0 \text{ and } P_\mu \left( \limsup_{t \uparrow \infty} e^{-\lambda_c t} Z_t(B) < \infty \right) = 1.$$

**Remark 2** (Total mass). In Theorem 3 we are concerned with the *local* behavior of the population size. When considering the *total mass* process $\|Z_t\| := Z_t(D)$, the growth rate may actually exceed $\lambda_c$. Indeed, take for example the $(L, \beta; D)$-branching diffusion with a conservative diffusion corresponding to $L$ on $D$ and with $\lambda_0 := \lambda_c(L, D) < 0$, and let $\beta > 0$ be constant. Then $\lambda_c(L + \beta, D) = \beta + \lambda_0 < \beta$, but – since the branching rate is spatially constant and since there is no 'loss of mass' by conservativeness – a classical theorem on Yule's processes tells us that $e^{-\beta t} \|Z_t\|$ tends to a nontrivial random variable as $t \to \infty$, that is, that the growth rate of the total mass is $\beta > \lambda_c$. ◇

Before the proofs we present a proposition, which is not hard to prove by an application of the Borel-Cantelli lemma.

**Proposition 2.** *For any nonempty open set $B \subset D$ and finite $\mu$,*

$$P_\mu \left( \limsup_{t \uparrow \infty} Z_t(B) \in \{0, \infty\} \right) = 1.$$

### 4.3. Proof of Theorem 3 (i)

**(a)** Assume that $\lambda_c \leq 0$. Then there exists an $h > 0$ solving $(L + \beta)h = 0$. We claim that $\{\langle h, Z_t \rangle : t \geq 0\}$ is a positive $P_\mu$-supermartingale for all $\text{supp}(\mu) \subset\subset D$. Indeed, it is not hard to show that $E_{\delta_x} \langle h, Z_t \rangle \leq h(x)$ for $t \geq 0$ and $x \in D$. From a standard application of the Markov property and the branching property, it then follows that $\{\langle h, Z_t \rangle : t \geq 0\}$ is a $P_\mu$-supermartingale.

Now in the possession of this supermartingale it follows for Borel $B \subset\subset D$, that

$$\limsup_{t \uparrow \infty} Z_t(B) \leq C \limsup_{t \uparrow \infty} \langle h, Z_t \rangle < \infty \tag{10}$$

$P_\mu$-almost surely where $C$ is a constant. When $B$ is open, by Proposition 2 it follows that $\lim_{t \uparrow \infty} Z_t(B) = 0$ $P_\mu$ – a.s. Since every compactly embedded Borel can be fattened up to an open $B \subset\subset D$, local extinction follows by comparison.

**(b)** Assume now that $\lambda_c > 0$. Since it is assumed $\text{supp}\,\mu \subset\subset D$, we can choose a large enough $B$ for which $\text{supp}\,\mu \subset B$ and $\lambda = \lambda_c(L + \beta, B) > 0$. Choose $0 \neq g \in C_c^+$ so that $g \leq \mathbf{1}_B$; obviously, it suffices to prove that $P_\mu \left( \limsup_{t \uparrow \infty} \langle g, Z_t \rangle > 0 \right) > 0$. Let $M_t^\phi$ and $\widetilde{P}_\mu$ be as in Theorem 2. By Lemma 1, it is enough to show that $\widetilde{P}_\mu \left( \limsup_{t \uparrow \infty} \langle g, Z_t \rangle > 0 \right) > 0$, or equivalently, that $\widetilde{E}_\mu \left( e^{-\limsup_{t \uparrow \infty} \langle g, Z_t \rangle} \right) < 1$. Let $\varepsilon > 0$. Use first Fatou's Lemma and then Theorem 2 to obtain the estimate

$$\widetilde{E}_\mu \left( e^{-\limsup_{t \uparrow \infty} \langle g, Z_t \rangle} \right) \leq \liminf_{t \uparrow \infty} \widetilde{E}_\mu \left( e^{-\langle g, Z_t \rangle} \right) \leq \liminf_{t \uparrow \infty} \mathbb{E}_{\phi\mu}^\phi \left( e^{-g(Y_t)} \right) \tag{11}$$

The ergodicity of $(Y, \mathbb{P}^\phi_{\phi\mu})$ implies that the right hand side of (11) is less than one, finishing the proof. □

Note the *intuition* behind the last part of the proof: the spine particle visits every part of $B$ for arbitrarily large times because it is an ergodic diffusion; this forces the process itself to do the same.

### 4.4. Proof of Theorem 3 (ii)

We may assume without loss of generality that $\lambda \in (0, \lambda_c)$. By standard theory, there exists a large enough $B^* \subset\subset D$ with a smooth boundary so that

$$\lambda^* := \lambda_c(L + \beta, B^*) \in (\lambda, \lambda_c].$$

Further, we can also choose $B^*$ large enough so that $\mathrm{supp}(\mu) \subset\subset B^*$.

We first claim that if $\Omega_0 := \{\lim_{t\uparrow\infty} e^{-\lambda t} Z_t(B^*) = \infty\}$, then $P_\mu(\Omega_0) > 0$. Indeed, if $\widehat{Z}$ is defined as the $(L, \beta; B^*)$-branching process then

$$P_\mu(\Omega_0) \geq P_\mu\left(\liminf_{t\uparrow\infty} e^{-\lambda^* t} Z_t(B^*) > 0\right) \geq P_\mu\left(\liminf_{t\uparrow\infty} e^{-\lambda^* t} \|\widehat{Z}_t\| > 0\right)$$
$$\geq P_\mu\left(\lim_{t\uparrow\infty} e^{-\lambda^* t} \langle \phi^*, \widehat{Z}_t \rangle > 0\right),$$

where $(L + \beta - \lambda^*)\phi^* = 0$ in $B^*$ and $\phi^* = 0$ on $\partial B^*$. Note it is implicit in the definition of $\widehat{Z}$ that particles are killed on the boundary $\partial B^*$. Since $\lambda^* > 0$, Lemma 1 implies that the last term in the last inequality is positive.

Now let $B$ be *any* open set with $\emptyset \neq B \subset\subset D$. Let $p := \inf_{x \in B^*} p(1, x, B) > 0$, where $\{p(t, ., dy) : t > 0\}$ is the transition measure for $(Y, \mathbb{P})$. Let $0 < q < p$ and $A_n := \{Z_{n+1}(B) \geq qZ_n(B^*)\}$. It follows from the law of large numbers and the Markov property that on $\Omega_0$, $\lim_{n\uparrow\infty} P(A_n \mid Z_n, ..., Z_1) = 1$. Using the extended Borel-Cantelli lemma, it easily follows that $\limsup_{t\uparrow\infty} e^{-\lambda t} Z_t(B) = \infty$ $P_{\delta_x}$−a.s. on $\Omega_0$. □

## 5. The Law of Large Numbers for spatial branching processes and superprocesses

In this section we would like to give a brief review on the problem of the Law of Large Numbers for spatial branching processes. Here $Z$ will denote the $(L, \beta, D)$-branching diffusion and $X$ will denote the $(L, \beta, \alpha, D)$-superprocess.

### 5.1. The question

Consider $Z$ and assume that it starts with a (finite) discrete measure $\mu$. Let us first assume that $Y$ is conservative (i.e. never reaches the cemetery state $\Delta$) and $\beta$ is a positive constant. Since the total mass process $|Z|$ is a non-spatial Yule's process (pure birth process), it is well known that $e^{-\beta t}|Z_t|$ has an a.s. positive limit as $t \to \infty$. That is, $|Z_t| \sim e^{\beta t} N_\mu$, where $N_\mu > 0$ a.s. Let $p^L(t, x, \cdot)$ denote the kernel corresponding to the diffusion $Y$ (or, equivalently, to $L$). It is natural to ask then whether it is true that for $B \subset\subset D$,
$$Z_t(B) \sim |Z_t| \cdot p^L(t, x, B)?$$
This would be a kind of Law of Large Numbers, because what it says is that the proportion of particles in $B$ at $t$ is given by the probability that a single particle is in $B$ at $t$. In other words, the question is whether
$$Z_t(B) \sim e^{\beta t} N_\mu \cdot p^L(t, x, B),$$
and since $e^{\beta t} \cdot p^L(t, x, \cdot) = p^{L+\beta}(t, x, \cdot)$ is the kernel corresponding to $L + \beta$, we can rewrite this in the form
$$Z_t(B) \sim N_\mu \cdot p^{L+\beta}(t, x, B).$$

Of course this last formula makes sense even for spatially dependent $\beta$'s. Moreover, it is known (and easy to derive) that $p^{L+\beta}(t, x, B)$ is just the expectation of $Z_t(B)$ starting from $\mu = \delta_x$ ('first moment formula'). So, for a general $\beta$ one has the following conjecture:
$$\frac{Z_t(B)}{E_{\delta_x} Z_t(B)} \sim N_x.$$
This conjecture is wrong! The reason is the possibility of local extinction. If local extinction occurs, then the left hand side becomes zero in finite time (and stays zero).

In terms of our spectral criterion for local extinction this means that there is no Law of Large Numbers, or more generally, there is no scaling for $Z_t(B)$ when $\lambda_c = \lambda_c(L + \beta) \leq 0$. The correct question therefore is this: is there a dichotomy in the sense that either
  (i) $Z$ suffers local extinction ($\lambda_c \leq 0$), or
  (ii) the (local) Law of Large Numbers holds for $Z (\lambda_c > 0)$?

This is not an easy question! In fact, since the spectral criterion of local extinction for superprocesses is the same, one can ask if the same dichotomy holds for the superprocess $X$. Moreover one has to clarify of course, what exactly '$\sim$' means in the formulas.

### *5.2. Some answers*

There are only partial answers in this direction. For superprocesses, it has been proven in [5] that if $D = \mathbb{R}^d$ and $\lambda_c > 0$ and furthermore a certain spectral condition holds for $L + \beta$, then for all $0 \neq g \in C_c^+$, the ratio $\frac{\langle X_t, g \rangle}{E_\mu \langle X_t, g \rangle}$ has a (non-degenerate) limit (depending only on $\mu$) in law.

One equivalent way of formulating the spectral condition is as follows. In general $p^{L+\beta}(t, x, B) \sim e^{\lambda_c t} \times$ possible subexponential term. The spectral condition is then equivalent to not having a subexponential term (the scaling is 'purely exponential'). For example, when $Y$ is $d$-dimensional Brownian motion and $\beta$ is positive constant, one has $p^{L+\beta}(t, x, B) \sim e^{\beta t} \cdot t^{-d/2}$, and therefore this case is not included in the setting. On the other hand, if one replaces Brownian motion with the Ornstein-Uhlenbeck process (or, in fact, with any positive recurrent process), then it is included.

It is important to point out however that, for nonconstant $\beta$'s, the underlying motion corresponding to $L$ *does not have to be positive recurrent* (it can even be transient), it is the operator $L + \beta$, involving branching too, that has to satisfy a certain spectral condition. In fact, an auxiliary diffusion, corresponding to a particular $h$-transform of the operator $L + \beta - \lambda_c$ is positive recurrent. This other diffusion also appears in a 'backbone'[12] construction.

This result was improved in [6] in the sense that a general Euclidean domain $D \subseteq \mathbb{R}^d$ was considered instead of $D = \mathbb{R}^d$ and convergence in law was upgraded to convergence in probability; recently it was further extended [4] to certain cases when the scaling is not purely exponential, including the case when $Y$ is $d$-dimensional Brownian motion and $\beta$ is positive constant.

Convergence in probability is far from being satisfying. In fact, by a well known theorem [1, 8], the analogous result for branching Brownian motion (with constant $\beta$) holds even in the strong sense (i.e. the convergence is a.s.), thus there is no real reason to assume that the convergence only holds in probability for superprocesses.

In three very recent papers [2, 3, 7] a.s. convergence was finally settled for certain classes of branching diffusions [3, 7] and superprocesses [2]. Despite the somewhat restrictive conditions, [2] seems to be the first time when SLLN is demonstrated for superprocesses, and it represents a major step in the direction of establishing the '*local extinction vs. SLLN*' dichotomy for superprocesses in terms of the associated generalized principal eigenvalue. Another nice feature of the paper is that besides symmetric diffusions, symmetric Levy processes are treated too. An essential difficulty is compounded in the fact that the method of [3] for replacing discrete times with continuous times does not carry through. The authors of [2] manage to overcome this difficulty with a highly nontrivial application of the martingale formulation for superprocesses.

At this point it seems to be a serious challenge to prove/ disprove the full dichotomy (without additional scaling assumptions) for both kinds of processes.

---

[12]'Backbone' constructions for superprocesses are similar to the 'spine' constructions for branching diffusions.

## 6. Further topics in superprocesses: compact support property and polar decomposition

### *6.1. The compact support property*

In this subsection we will consider superprocesses with a somewhat more general branching mechanism. For simplicity, we take $D = \mathbb{R}^d$, although most of the results are true for general Euclidean domains. [13]

In this case the evolution equation appearing in the log-Laplace formulation is of the form

$$\begin{aligned} u_t &= Lu + \beta u - \alpha u^p \ \text{ in } \mathbb{R}^d \times (0, \infty); \\ u(x, 0) &= f(x) \ \text{ in } \mathbb{R}^d; \\ u(x, t) &\geq 0 \ \text{ in } \mathbb{R}^d \times [0, \infty), \end{aligned} \quad (12)$$

where $p \in (1, 2]$. When $1 < p < 2$ the variance is no longer finite and only the moments up to $p$ exist in the particle process approximation, whereas the scaling of the clock becomes $n^{p-1}$. Finally, the paths of the superprocess are not continuous but only càdlàg when $1 < p < 2$.

Let $\mu \in \mathcal{M}_c(\mathbb{R})$. Since $D = \mathbb{R}^d$, the measure-valued process corresponding to $P_\mu$ possesses the *compact support property* if

$$P_\mu \left( \bigcup_{0 \leq s \leq t} \text{supp } X(s) \text{ is bounded} \right) = 1, \text{ for all } t \geq 0. \quad (13)$$

The definition is independent of the choice of the nonzero compactly supported starting measure.

There are four objects, corresponding to four different underlying probabilistic effects, which can influence the compact support property:

1. $L$, the operator corresponding to the underlying motion;
2. $\beta$, the mass creation parameter of the branching mechanism;
3. $\alpha$, the nonlinear component of the branching mechanism, which can be thought of as the variance parameter if $p = 2$;
4. $p$, the power of the nonlinearity, which is the scaling power and is connected to the fractional moments of the offspring distribution in the particle process approximation to the measure-valued process.

We shall see that both $L$ and $\alpha$ play a large role in determining whether or not the compact support property holds; $\beta$ and $p$ play only a minor role.

The connection with partial differential equations is that the compact support property turns out to be equivalent to a uniqueness property for solutions to (12).

**Result 1.** *The compact support property holds for one, or equivalently all, nonzero, compactly supported initial measures $\mu$ if and only if there are no nontrivial solutions to (12) with initial data $f \equiv 0$.*

---

[13]The results of this subsection are joint with R. Pinsky.

**Remark 3.** Result 1 suggests a parallel between the compact support property for measure-valued processes and the non-explosion property for diffusion processes. Indeed, the non-explosion property for the diffusion process $Y$ corresponding to the operator $L$ is equivalent to the nonexistence of nontrivial, *bounded* positive solutions to the linear Cauchy problem with 0 initial data:

$$u_t = Lu \text{ in } \mathbb{R}^d \times (0, \infty);$$
$$u(x, 0) = 0 \text{ in } \mathbb{R}^d; \qquad (14)$$
$$u \geq 0 \text{ in } \mathbb{R}^d \times (0, \infty).$$

It is natural for bounded, positive solutions to be the relevant class of solutions in the linear case and for positive solutions to be the relevant class of solutions in the semi-linear case. Indeed, by Ito's formula, the probabilities for certain events related to $Y$ are obtained as bounded, positive solutions to the linear equation, and by the log-Laplace equation, the negative of the logarithm of the probability of certain events related to $X$ can be obtained as positive solutions to the semi-linear equation. ◇

The class of operators $L$ satisfying the following assumption will play an important role.

**Assumption 1.** *For some $C > 0$,*

1. $\sum_{i,j=1}^n a_{ij}(x)\nu_i\nu_j \leq C|\nu|^2(1+|x|^2), \ x, \nu \in \mathbb{R}^d;$
2. $|b(x)| \leq C(1+|x|), \ x \in \mathbb{R}^d.$

We have the following basic result.

**Result 2.** *Let $p \in (1, 2]$ and let the coefficients of $L$ satisfy Assumption 1.*

1. *There is no nontrivial solution to (14); thus, the diffusion process $Y$ does not explode.*
2. *If*
$$\inf_{x \in \mathbb{R}^d} \alpha(x) > 0,$$
   *then there is no nontrivial solution to (12) with initial data $f = 0$; thus, the compact support property holds for $X$.*

The conditions in Assumption 1 are classical conditions which arise frequently in the theory of diffusion processes. Result 2 shows that if the coefficients of $L$ obey this condition and if the branching coefficient $\alpha$ is bounded away from zero, then everything is well behaved—neither can the underlying diffusion process explode nor can the measure-valued process fail to possess the compact support property.

The following result shows that the compact support property can fail if $\inf_{x \in \mathbb{R}^d} \alpha(x) = 0$. It also demonstrates that the effect of $\alpha$ on the compact support property cannot be studied in isolation, but in fact depends on the underlying diffusion.

**Theorem 4.** *Let $p \in (1, 2]$ and let*

$$L = \frac{1}{2} \sum_{i,j=1}^{d} a_{i,j} \frac{\partial^2}{\partial x_i \partial x_j}, \quad \text{where}$$

$$C_0^{-1}(1 + |x|)^m \leq \sum_{i,j=1}^{d} a_{i,j}(x) \leq C_0(1 + |x|)^m, \; m \in [0, 2], \; \text{for some } C_0 > 0.$$

1. *If*
$$\alpha(x) \geq C_1 \exp(-C_2|x|^{2-m}),$$
   *for some $C_1, C_2 > 0$, then the compact support property holds for $X$.*
2. *If*
$$\alpha(x) \leq C \exp(-|x|^{2-m+\epsilon}) \text{ and } \beta(x) \geq -C(1+|x|)^{2-m+2\delta},$$
   *for some $C, \epsilon > 0$ and some $\delta < \epsilon$, then the compact support property does not hold for $X$.*

(By Result 1, to prove this theorem, it is necessary and sufficient to show that if $\alpha$ is as in part (1) of the theorem, then there is no nontrivial solution to (12) with initial data $f = 0$, while if $\alpha$ is as in part (2) of the theorem, then there is such a nontrivial solution.)

As a complement to Theorem 4, we note the following result.

**Result 3.** *Let $p \in (1, 2]$.*

1. *Let $d \geq 2$ and let*
$$L = A(x)\Delta, \; \text{where } A(x) \geq C(1 + |x|)^m, \; \text{for some } C > 0 \text{ and } m > 2.$$
   *Assume that*
   $$\sup_{x \in \mathbb{R}^d} \alpha(x) < \infty \text{ and } \beta \geq 0.$$
   *Then the compact support property does not hold for $X$.*
2. *Let $d = 1$ and let*
$$L = A(x)\frac{d^2}{dx^2}, \; \text{where } A(x) \geq C(1+|x|)^m, \; \text{for some } C > 0 \text{ and } m > 1+p.$$
   *Assume that*
   $$\sup_{x \in \mathbb{R}^d} \alpha(x) < \infty \text{ and } \beta \geq 0.$$
   *Then the compact support property does not hold for $X$.*
3. *Let $d = 1$ and let*
$$L = A(x)\frac{d^2}{dx^2}, \; \text{where } A(x) \leq C(1+|x|)^m, \; \text{for some } C > 0, \text{ and } m \leq 1+p.$$
   *Assume that*
   $$\inf_{x \in \mathbb{R}^d} \alpha(x) > 0 \text{ and } \beta \leq 0.$$
   *Then the compact support property holds for $X$.*

**Remark 4.** It follows from Result 3 that if $d = 1$ and $L = (1 + |x|)^m \frac{d^2}{dx^2}$, with $m \in (2, 3]$, and say $\alpha = 1$ and $\beta = 0$, then the compact support property will depend on the particular choice of $p \in (1, 2]$.      ⋄

**Remark 5.** If $d = 2$ and $L = (1 + |x|)^m \Delta$, with $m > 2$, and say $\alpha = 1$ and $\beta = 0$, then by Result 3, the compact support property does not hold, yet the underlying diffusion does not explode since it is a time-change of a recurrent process; namely, of two-dimensional Brownian motion.      ⋄

The proof is based on the fact that, by Result 1, it is necessary and sufficient to show that under the conditions of parts (1) and (2), there exists a nontrivial solution to (12), while under the conditions of part (3) there does not.

A useful comparison result is as follows.

**Proposition 3.** *Assume that*
$$\beta_1 \leq \beta_2$$
*and*
$$0 < \alpha_2 \leq \alpha_1.$$
*If the compact support property holds for $\beta = \beta_2$ and $\alpha = \alpha_2$, then it also holds for $\beta = \beta_1$ and $\alpha = \alpha_1$.*

Theorem 4 and Result 3 demonstrate the effect of the underlying diffusion $Y$ on the compact support property in the case that $L$ is comparable to $(1+|x|)^m \Delta$. We now consider more generally the effect of the underlying diffusion process on the compact support property.

**Theorem 5.** *Let $p \in (1, 2]$ and assume that the underlying diffusion process $Y$ explodes. Assume that*
$$\sup_{x \in \mathbb{R}^d} \alpha(x) < \infty \text{ and } \inf_{x \in \mathbb{R}^d} \beta(x) > -\infty. \tag{15}$$
*Then the compact support property does not hold.*

**Remark 6.** Theorem 5 shows that if the branching mechanism satisfies (15), then the compact support property never holds if the underlying diffusion is explosive. The converse is not true—an example was given in Remark 2 following Result 3, and another one appears in the remark following Corollary 3.     ⋄

The next result shows that the restriction $\sup_{x \in \mathbb{R}^d} \alpha(x) < \infty$ in Theorem 5 is essential.

**Proposition 4.** *Let $p \in (1, 2]$. Let $m \in (-\infty, \infty)$,*
$$L = (1 + |x|)^m \Delta \text{ in } \mathbb{R}^d,$$
*$\beta = 0$ and*
$$\alpha(x) \geq c(1 + |x|)^{m-2},$$
*for some $c > 0$. Then the compact support property holds for the measure-valued process $X$. However, if $m > 2$ and $d \geq 3$, the diffusion process $Y$ explodes.*

We now discuss how the concept of a measure-valued process *hitting* a point can be formulated and understood in terms of the compact support property.

Similarly to the range define also $\mathcal{R}_t = \overline{\left(\cup_{s\in[0,t]} \text{supp}(X(s))\right)}$. A path of the measure-valued process is said to *hit* a point $x_0 \in \mathbb{R}^d$ if $x_0 \in \mathcal{R}$. If $X$ exhibits local extinction, then $x_0 \in \mathcal{R}$ if and only if $x_0 \in \mathcal{R}_t$ for sufficiently large $t$. Thus, we have:

If $X$ exhibits local extinction, then

$$P_\mu(X \text{ hits } x_0) > 0 \text{ if and only if there exists a } t > 0 \text{ such that} \quad (16)$$
$$P_\mu(\cup_{0\leq s\leq t} \text{supp}(X(s)) \text{ is not compactly embedded in } \mathbb{R}^d - \{x_0\}) > 0.$$

Now although we have assumed in this section that the underlying state space is $\mathbb{R}^d$, everything goes through just as well on an arbitrary domain $D \subset \mathbb{R}^d$. Recall that the compact support property is defined with respect to the domain $D$, and the underlying diffusion will explode if it hits $\Delta$ in finite time. In particular, Result 1 still holds with $\mathbb{R}^d$ replaced by $D$.

In light of the above observations, consider a measure-valued process $X$ corresponding to the log-Laplace equation (12) on $\mathbb{R}^d$ with $d \geq 2$. The underlying diffusion process $Y$ on $\mathbb{R}^d$ corresponds to the operator $L$ on $\mathbb{R}^d$. Let $\hat{Y}$ denote the diffusion process on the domain $D = \mathbb{R}^d - \{x_0\}$ with absorption at $x_0$ and corresponding to the same operator $L$. (Note that if $x_0$ is polar for $Y$, then $Y$ and $\hat{Y}$ coincide when started from $x \neq x_0$. In fact, $x_0$ is always polar under the assumptions we have placed on the coefficients of $L$.) Let $\hat{X}$ denote the measure-valued process corresponding to the log-Laplace equation (12), but with $\mathbb{R}^d$ replaced by $D = \mathbb{R}^d - \{x_0\}$. *It follows from (16) that if $X$ suffers local extinction, then the measure-valued process $X$ hits the point $x_0$ with positive probability if and only if $\hat{X}$ on $\mathbb{R}^d - \{x_0\}$ does not possess the compact support property.* Furthermore, the above discussion shows that even if $X$ does not suffer local extinction with probability one, a sufficient condition for $X$ to hit the point $x_0$ with positive probability is that $\hat{X}$ on $\mathbb{R}^d - \{x_0\}$ does not possess the compact support property.

A similar analysis can be made when $d = 1$. (The process $\hat{X}$ above must be replaced by two processes, $\hat{X}^+$ and $\hat{X}^-$, defined respectively on $(x_0, \infty)$ and $(-\infty, x_0)$.) In the sequel, we will assume that $d \geq 2$.

Consider now the following semi-linear equation in $\mathbb{R}^d - \{0\}$:

$$\begin{aligned} u_t &= \frac{1}{2}\Delta u - u^p \text{ in } (\mathbb{R}^d - \{0\}) \times (0, \infty); \\ u(x,0) &= 0 \text{ in } \mathbb{R}^d - \{0\}; \\ u &\geq 0 \text{ in } (\mathbb{R}^d - \{0\}) \times [0, \infty). \end{aligned} \quad (17)$$

By Result 1, the measure-valued process corresponding to the semi-linear equation $u_t = \frac{1}{2}\Delta u - u^p$ in $\mathbb{R}^d - \{0\}$ will possess the compact support property if and only if (17) has a nontrivial solution. The following result was recently proved by P. Baras, and M. Pierre [1], and later by R. Pinsky [6], using a different method.

**Result 4.** *Let $p > 1$ and $d \geq 2$.*

1. *If $d < \frac{2p}{p-1}$, then there exists a nontrivial solution to (17).*
2. *If $d \geq \frac{2p}{p-1}$, then there is no nontrivial solution to (17).*

Here $p > 1$ is unrestricted, even though of course there is no probabilistic import when $p > 2$.

In terms of the superprocess, the following was shown in [4] (see also the remark at the end of this subsection).

**Theorem 6.** *Let $p > 1$ and $d \geq 2$. Let $X$ denote the measure-valued process on the* punctured *space $\mathbb{R}^d - \{0\}$ corresponding to the semi-linear equation $u_t = \frac{1}{2}\Delta u + \beta u - u^p$ in $(\mathbb{R}^d - \{0\}) \times (0, \infty)$. Consider the Cauchy problem*

$$\begin{aligned} u_t &= \frac{1}{2}\Delta u + \beta u - u^p \text{ in } (\mathbb{R}^d - \{0\}) \times (0, \infty); \\ u(x, 0) &= 0 \text{ in } \mathbb{R}^d - \{0\}; \\ u &\geq 0 \text{ in } (\mathbb{R}^d - \{0\}) \times [0, \infty). \end{aligned} \tag{18}$$

*Assume that*
$$d < \frac{2p}{p-1}.$$

*Let*
$$\beta_0 = \frac{d(p-1) - 2p}{(p-1)^2} < 0.$$

1. *If*
$$\beta(x) \geq \frac{\beta_0 + \kappa}{|x|^2}, \text{ for some } \kappa \in (0, -\beta_0],$$

   *then there exists a nontrivial solution to (18); hence, the compact support property does not hold for $X$.*
2. *If*
$$\limsup_{x \to 0} |x|^2 \beta(x) < \beta_0,$$

   *then there is no nontrivial solution to (18); hence, the compact support property holds for $X$.*

**Remark 7.** The restriction $\kappa \leq -\beta_0$ is made to ensure that $\beta$ is bounded from above. ◇

We have the following corollary of Result 4 and Theorem 6.

**Corollary 3.** *Let $X$ denote the measure-valued process on all of $\mathbb{R}^d$, $d \geq 2$, corresponding to the semi-linear equation $u_t = \frac{1}{2}\Delta u + \beta u - u^p$ in $\mathbb{R}^d \times (0, \infty)$.*

1. *If $\beta$ is bounded from below and $d < \frac{2p}{p-1}$, then $X$ hits any point $x_0$ with positive probability;*
2. *If $\beta \leq 0$ and $d \geq \frac{2p}{p-1}$, then $X$ hits any point $x_0$ with probability 0;*

3. If $\beta \leq 0$, $d < \frac{2p}{p-1}$ and $\beta$ has a singularity at the origin such that

$$\limsup_{x \to 0} |x|^2 \beta(x) < \frac{d(p-1) - 2p}{(p-1)^2},$$

then $X$ hits 0 with probability 0.

When $\beta = 0$, part (1) follows immediately from Result 4, Result 1 and the discussion preceding Result 4. The general case requires a little more works as one needs to compare different $\beta$'s. (Recall that we are always assuming here that $\beta$ is bounded.)

When $\beta \leq 0$, it is well-known that the super-Brownian motion in the statement of the corollary becomes extinct. Thus, part (2) follows from Result 4, Result 1 and the discussion preceding Result 4, while part (3) follows from Theorem 6, Result 1 and the discussion preceding Result 4.

**Remark 8.** When $\beta = 0$, the results in parts (1) and (2) of Corollary 3 state that critical, super-Brownian motion hits a point with positive probability if $d < \frac{2p}{p-1}$, and with zero probability if $d \geq \frac{2p}{p-1}$. (Note that this also yields an example where $X$ does not possess the compact support property even though the underlying diffusion process does not explode.[14]) This result was first proved by D. Dawson, I. Iscoe, and E. Perkins [2], and by E. B. Dynkin and S. Kuznetsov [3]. Their methods were different than the one presented here. ◇

## *6.2. Polar decomposition*

As we have already seen in the beginning of the previous section, the total mass process of a branching diffusion is sometimes easier to handle than the whole process. This is true for superprocesses too. In fact, one can even think of the superprocess until its extinction as a process of two components: the total mass process, as 'radial component' and the process normalized by the total mass (which is a probability measure valued process) as the 'angular component'. We are going to demonstrate this for the critical super-Brownian motion (the $(\frac{1}{2}\Delta, 0, 1, R^d)$-superprocess).

To this end, let us have a review first on the *Fleming-Viot superprocess* which is of independent interest too. Imagine that we have $N$ individuals on the lattice $\mathbb{Z}^d$ and each of them perform independent symmetric random walk. In the genetic setting $\mathbb{Z}^d$ is the space of genetic types and the random walks are *mutations*. Each individual has her own independent exponential clock with rate $\gamma > 0$ (*sampling rate*) and when the clock rings, she jumps to a position randomly chosen from the current empirical distribution of the population (*sampling mechanism*). Note that we can consider this mechanism as branching: one individual dies and another one is born somewhere. Also, the process can be viewed as a (discrete) measure valued process.

---

[14]It is easy to find the connection and we leave this to the reader.

Let us *rescale* now the model as follows: assign mass $1/N$ to each individual (and so the measure valued process becomes a probability measure valued process). Furthermore speed up time by factor $N$ and reduce lattice mesh by factor $1/\sqrt{N}$. (So the mutation part will converge to Brownian motion).

Formally, let $V_N(t, A) := \sum_{j/\sqrt{N} \in A} p(Nt, j)$, where $A \subset \mathbb{Z}^d$ and $p(s, j)$ is the proportion of the population at the point $j$ at time $s \geq 0$. Let $N \to \infty$ now. If $V_N(0, \cdot)$ converges to the probability measure $\nu$ on $R^d$ weakly, then the process $\{V_N(t, \cdot); t \geq 0\}$ converges to a process $\{V(t, \cdot); t \geq 0\}$ weakly. This limiting (probability measure valued) process is the Fleming-Viot superprocess.

Let us write $X$ in a 'polar coordinates': $X_t = r_t \cdot \theta_t$, where $r_t := |X_t|$ and $\theta_t := \frac{X_t}{|X_t|}$ for $r_t > 0$, (otherwise $r_t = 0, X_t = \mathbf{0}$). Then it is easy to show that $r$ is a diffusion process ('Feller's diffusion') on $D = [0, \infty)$ corresponding to the operator $L = x \frac{d^2}{dx^2}$. (This diffusion will hit zero with probability one and afterwards it stays at zero.) A result due to P. March and A. Etheridge then states that $\theta$ is a Fleming-Viot superprocess, except that the sampling rate is replaced by $1/r_t$. This 'skew-product' representation is the polar decomposition for critical super-Brownian motion.

## 7. Classical problems for random media: Brownian motion among Poissonian obstacles and related topics

In this section we give a brief (and surely not comprehensive) review on some of the classical problems in random media.

### 7.1. Brownian survival among Poissonian traps

As usual, $B(x, r)$ denotes the open ball of radius $r$ around $x$. Let $W$ denote the $d$-dimensional Brownian motion (Wiener-process) starting at the origin (the corresponding probability will be denoted by P) and let us define the *Wiener-sausage* up to time $t$ with radius $a$:

$$V_t^a = \text{Vol}\left(\bigcup_{0 \leq s \leq t} \overline{B}(W_t, a)\right).$$

It turns out that this object is intimately related to the problem of *Brownian survival among Poissonian traps*. Let $\omega$ denote the Poisson point process (PPP) with intensity $\nu(dx)$ on $\mathbb{R}^d$. (The corresponding probability will be denoted by $\mathbb{P}$.) That is, the number of points in $A \subset \mathbb{R}^d$ bounded is Poisson with parameter $\nu(A)$ and the number of points in $A \subset \mathbb{R}^d$ and the number of points in $B \subset \mathbb{R}^d$ are independent random variables if $A$ and $B$ are disjoint. So for example $\mathbb{P}(A \text{ is empty}) = e^{-\nu(A)}$. Unless indicated otherwise we will set $\nu(dx) = \nu \cdot dx$, $\nu > 0$.

The (random) set
$$K := \bigcup_{x_i \in \text{supp}(\omega)} \overline{B}(x_i, a)$$
is a *trap configuration*, or *obstacle*.

If $T$ denotes the first hitting time of $K$ by our Brownian motion, then we would like to know what the asymptotic behavior of the probability of $\{T > t\}$ is. Of course this question can be asked in two different ways: what the asymptotics of $(\mathbb{E} \otimes P)(T > t)$ is, or what the $\mathbb{P}$-a.s. asymptotics of $P(T > t)$ is. In the first case one talks about the *annealed* asymptotics and in the second one, about the *quenched* one.

Now, the annealed asymptotics is in fact just the same as the asymptotics for the Laplace transform of the Wiener sausage. Indeed,

$$(\mathbb{E} \otimes P)(T > t) = (\mathbb{E} \otimes P)\left(\bigcup_{0 \leq s \leq t} \overline{B}(W_t, a) \text{ is empty}\right) = E e^{-\nu V_t^a}.$$

A little more than 30 years ago Donsker and Varadhan proved that for all $a, \nu > 0$,
$$\lim_{t \to \infty} t^{-d/(d+2)} \log E \exp(-\nu V_t^a) = -\tilde{c}(d, \nu),$$

where $\tilde{c}(d,\nu)$ is an explicitly given positive constant (we are using $\tilde{c}$ because there is another constant $c(d,\nu)$ in the quenched case). Equivalently,

$$Ee^{-\nu V_t^a} = \exp[-\tilde{c}(d,\nu)t^{d/(d+2)}(1+o(1))].$$

Note that the scaling is milder than exponential, it is often referred to as '*stretched exponential*'. For example in two dimension the Laplace transform behaves like $e^{-\text{const}\sqrt{t}}$.

Of course, the same is then true for the annealed survival probability:

$$(\mathbb{E} \otimes P)(T > t) = \exp[-\tilde{c}(d,\nu)t^{d/(d+2)}(1+o(1))].$$

The easy part is

$$\liminf_{t\to\infty} t^{-d/(d+2)} \log E\exp(-\nu V_t^a) \geq -\tilde{c}(d,\nu),$$

because all one has to do is to 'advise' a survival strategy which is realized (at least) with probability $\exp[-\tilde{c}(d,\nu)t^{d/(d+2)}(1+o(1))]$, whereas in order to verify

$$\limsup_{t\to\infty} t^{-d/(d+2)} \log E\exp(-\nu V_t^a) \leq -\tilde{c}(d,\nu),$$

one has to show that the advised strategy is actually optimal.

In the context of Brownian survival among Poissonian traps these (and also the quenched asymptotics) were settled by A. S. Sznitman.

Let us see now how to 'advise' a survival strategy. An obvious observation is that if $B(0, R_t + a)$ is empty (such an empty ball is called a *clearing*), and if $W$ is confined to $B(0, R_t)$ up to $t$, then $T > t$. That is,

$$(\mathbb{E} \otimes P)(T > t) \geq \mathbb{P}(B(0, R_t + a) \text{ is empty}) \cdot P(\tau_{B(0,R_t)} > t),$$

where $\tau_{B(0,R_t)}$ denotes the exit time of $W$ from the ball. It is not hard to show that $P(\tau_{B(0,R_t)} > t) \geq c(d)e^{-\lambda_d/R_t^2 \cdot t}$, where $\lambda_d$ is the principal Dirichlet eigenvalue of $-\Delta$ on the unit ball. Hence

$$(\mathbb{E} \otimes P)(T > t) \geq c(d)e^{-\nu\omega_d(R_t+a)^d - \lambda_d/R_t^2 \cdot t},$$

where $\omega_d$ is the volume of the $d$-dimensional unit ball. So

$$\log(\mathbb{E} \otimes P)(T > t) \geq \log c(d) - \nu\omega_d(R_t + a)^d - \lambda_d/R_t^2 \cdot t. \tag{19}$$

Of course, one wants to get the optimal (strongest) lower estimate. In order to achieve this, we match the orders of the last two terms (the constant will be lower order and thus negligible), which happens when $R_t = \mathcal{O}(t^{1/(d+2)})$. If $R_t = kt^{1/(d+2)}$, then we have to maximize $f(k) := \nu\omega_d k^d + \frac{\lambda_d}{k^2}$, which is found easily (by differentiation): $k = \left(\frac{2\lambda_d}{d\nu\omega_d}\right)^{1/(d+2)}$. The rest is straightforward. Setting $R_t := \left(\frac{2\lambda_d}{d\nu\omega_d}\right)^{1/(d+2)} \cdot t^{1/(d+2)}$ in (19), dividing by $t^{d/(d+2)}$ and finally taking liminf gives the desired lower estimate.

As we have mentioned, the proof of the upper estimate is much harder and it is definitely beyond the scope of these notes.

As far as the quenched problem is concerned, the asymptotics is different. If $a$ is sufficiently small such that the origin belongs to an infinite trap free cluster with positive $\mathbb{P}$-probability, then on that set $\mathbb{P}$-a.s,

$$P(T > t) = \exp\left(-c(d,\nu)\frac{t}{(\log t)^{2/d}}(1 + o(1))\right),$$

where again $c(d,\nu)$ is an explicitly given positive constant.

Finally let us note that the annealed and quenched asymptotics hold also for *soft* obstacles (so far we have been talking about *hard* ones), where the meaning of being 'soft' is as follows. Let $b > 0$ (the size of $b$ plays no role) and for a given realization of $\omega$ let

$$P(T > t) = \exp\left(-b\int_0^t \mathbb{1}_K(W_s)\, ds\right).$$

(This is not the usual formulation because we do not add up the rates on overlapping balls, but it does not matter from the point of view of asymptotics and this is what we will need.) That is, as long as $W$ is inside $K$ it 'feels' an exponential killing rate, whereas in the 'Swiss cheese' $K^c$ it is 'safe'.

## 7.2. Random walk in random environment

Let us first assume $d = 1$. There are two basic models for random walk in random environment (RWRE):

- **Site randomness** (this is what usually meant by RWRE). For a given realization of the environment $\omega$, $\{X_n\}$ is a Markov-chain on $\mathbb{Z}$. The walker located at $x \in \mathbb{Z}$ moves to $x+1$ with probability $p(x,\omega)$ and to $x-1$ with probability $q(x,\omega) = 1 - p(x,\omega)$, where $\{p(x,\omega),\ x \in \mathbb{Z}\}$ are i.i.d. random variables with values in $[0,1]$.
- **Bond randomness** (sometimes called 'random conductivity' model). Here $\{c_{x,x+1}(\omega),\ x \in \mathbb{Z}\}$ are i.i.d. with values in $(0,\infty)$ and for a given realization $\{X_n\}$ is a Markov-chain on $\mathbb{Z}$ with probability

$$p(x,\omega) := \frac{c_{x,x+1}(\omega)}{c_{x-1,x}(\omega) + c_{x,x+1}(\omega)}$$

of jumping to $x+1$ and with probability $q(x,\omega) = 1 - p(x,\omega)$ of jumping to $x-1$, in which case the $\{p(x,\omega),\ x \in \mathbb{Z}\}$ are of course no longer i.i.d..

Site randomness goes back to Chernov [2] in the early sixties, while bond randomness goes back to Fatt [3] in the fifties.

Of course it is no problem to define these models in higher dimensions too. However, the higher dimensional case is significantly harder than the one-dimensional one, and much less is known. In the sequel we keep assuming $d = 1$.

One of the natural questions concern the speed of the random walk. Now, in the bond randomness model one can show that for almost all realizations of $\omega$, $\frac{X_n}{n}$ tends to zero, which is interpreted as 'vanishing limiting velocity'. However, for the site randomness case, the situation is far more subtle. Here we suffice with citing the classical result of Solomon [4].

**Theorem 7.** *Let $\rho(x,\omega) := \frac{q(x,\omega)}{p(x,\omega)}$ and $\rho := \rho(0,\omega)$. If $\mathbb{P}$ corresponds to the environment and $P^\omega$ to the walk given $\omega$, then there are three cases.*

1. $\mathbb{E}[\log \rho] < 0$ and $\lim_n X_n = +\infty$,
2. $\mathbb{E}[\log \rho] > 0$ and $\lim_n X_n = -\infty$,
3. $\mathbb{E}[\log \rho] = 0$ and $\liminf_n X_n = -\infty$ but $\limsup_n X_n = +\infty$,

*where all the limits are meant $\mathbb{P} \otimes P^\omega$-a.s.*

*Furthermore, $\frac{X_n}{n} \to v$ $\mathbb{P} \otimes P^\omega$-a.s. as $n \to \infty$, where, again, there are three cases.*

1. $\mathbb{E}\rho < 1$ and $v = \frac{1-\mathbb{E}\rho}{1+\mathbb{E}\rho} \in (0,1)$,
2. $\mathbb{E}\rho \geq 1$, $\mathbb{E}[\rho^{-1}] \geq 1$ and $v = 0$,
3. $\mathbb{E}\rho^{-1} < 1$ and $v = \frac{\mathbb{E}\rho^{-1}-1}{1+\mathbb{E}\rho^{-1}} \in (-1,0)$.

**Remark 9.** To understand the statements better, note the following.

**(i)** In the third case $\mathbb{E}\rho^{-1} < 1$ automatically implies $\mathbb{E}\rho > 1$, by Jensen. Also by Jensen, $\mathbb{E}\rho < 1$ implies $\mathbb{E}[\log \rho] < 0$ and $\mathbb{E}\rho^{-1} < 1$ implies $\mathbb{E}[\log \rho] > 0$.

**(ii)** When $v = 0$, it is known that the propagation can be as slow as $(\log n)^2$ per $n$ steps.

**(iii)** It is easy to show (again by Jensen), that $\mathbb{E}\rho < 1$ implies $v \leq \mathbb{E}p - \mathbb{E}q$, where $p = p(0,\omega)$ and $q = 1-p$. Similarly $\mathbb{E}\rho^{-1} < 1$ implies $v \geq \mathbb{E}p - \mathbb{E}q$. These inequalities are strict (unless $p$ is deterministic), which shows the 'slowdown' of the walk relative to the naive guess that the speed is simply $v^* := \mathbb{E}p - \mathbb{E}q$.  ⋄

## 8. Spatial branching processes among Poissonian obstacles: hard obstacles

We now start putting Poissonian obstacles and spatial branching together.[15]

Let $T$ denote the hitting time of $K$ (defined in the previous section) by $Z$, where $Z$ denotes branching Brownian motion (BBM) with strictly dyadic branching and constant rate $\beta > 0$, starting with a single particle at the origin. That is $T$ is the time the first particle hits any trap:

$$T := \inf\{s \geq 0, \ W_s \in K\}.$$

Our question is again the asymptotic behavior of the probability of $\{T > t\}$.

### 8.1. Constant trap intensity

Let $\mathbb{P}$ be as in the previous section and let $P$ denote the probability corresponding to the BBM.

**Theorem 8.** *When $d \geq 2$, the* annealed *asymptotics is as follows:*

$$\lim_{t \to \infty} \frac{1}{t} \log(\mathbb{P} \otimes P)(T > t) = -\beta,$$

*and this is also true for the supercritical super-Brownian motion starting with unit mass at the origin. More precisely, if $X$ is the $(\frac{1}{2}\Delta, \beta, \alpha; \mathbb{R}^d)$-superprocess with $\alpha, \beta$ positive constants and $d \geq 2$, then*

$$\lim_{t \to \infty} \frac{1}{t} \log(\mathbb{P} \otimes P_{\delta_0})(T^* > t \mid S) = -\beta,$$

*where $P$ corresponds to $X$, $T^* := \inf\{t \geq 0 \mid X_t(K) > 0\}$ and $S$ is the event that $X$ survives forever.*

Although this result is far from being trivial, in a sense it is not too exciting because it indicates the overwhelming influence of the branching over the spatial motion. For example, in case of $Z$, the probability of not branching at all up to $t$ is $e^{-\beta t}$ and by the Donsker-Varadhan result we know that a single particle avoids $K$ up to $t$ with a subexponentially small probability. Therefore, we immediately obtain the lower bound

$$\liminf_{t \to \infty} \frac{1}{t} \log(\mathbb{P} \otimes P)(T > t) \geq -\beta,$$

for all $d \geq 1$. That this actually coincides with the upper bound for $d \geq 2$, shows that in a sense suppressing the branching is the 'best the system can do' in order to avoid the traps, and the cost of doing so is the dominating factor on a logarithmic scale. Because of the conditioning for $X$ however, the lower estimate for the superprocess is not easy.

---

[15] Results of the first subsection are due to the author, results in the second one are joint with F. den Hollander.

### 8.2. Decaying trap intensity

In order to treat the more delicate one-dimensional case of the previous subsection, we are going to discuss a setting when the trap intensity decays at the order $|x|^{d-1}$ as $|x| \to \infty$, which includes are previous one-dimensional case. So let us consider now a trap intensity satisfying that

$$\frac{\mathrm{d}\nu}{\mathrm{d}x} \sim \ell |x|^{1-d},$$

where $\ell$ is a fine tuning constant. Considering the same process $Z$ under $P$ as before and the same killing rule (we kill $Z$ when the first particle hit $K$) we have the following annealed asymptotics for the tail of the survival probability.

**Theorem 9.**
$$\lim_{t\to\infty} \log \frac{1}{t} (\mathbb{P} \otimes P)(T > t) = -I(\ell, \beta, d),$$

where $T$ is as before and $I(\ell, \beta, d) > 0$ can be computed in terms of a variational problem.

Instead of writing down the variational problem, let us see the *optimal survival strategy* the variational problem 'encodes'.

First, it turns out that there is a *crossover* at a critical value of $\ell$, which we will denote by $\ell_{cr} = \ell_{cr}(\beta, d)$. (This crossover is the reason we chose this particular order for $\nu$. If the density is higher or lower order than $|x|^{1-d}$, then we fall into one of the regimes already present when $\frac{\mathrm{d}\nu}{\mathrm{d}x} \sim \ell|x|^{1-d}$, and $\ell > \ell_{cr}$ or $\ell < \ell_{cr}$, respectively.) If $\ell_{cr}^* := s_d^{-1}\sqrt{\beta/2}$, where $s_d$ is the surface of the $d$-dimensional unit ball, then for $d = 1$, $\ell_{cr} = \ell_{cr}^*$, while for $d \geq 2$, $\ell_{cr} < \ell_{cr}^*$. In the theorem below we try to explain the optimal survival strategy in words rather than by using the precise formulation. Roughly speaking, the optimal survival strategy means that, conditioned on survival, the system, consisting of the BBM as well as the traps, behaves in a certain way up to time $t$ with $\mathbb{P} \otimes P$-probability tending to one as $t \to \infty$. We will distinguish between a 'low intensity' ($\ell < \ell_{cr}$) and a 'high intensity' ($\ell > \ell_{cr}$) regime. Note, that the case $\ell = \ell_{cr}$ is left open.

Concerning the optimal survival strategy, we have the following result. (The first part is actually just an interpretation of the precise statement which will be given in the theorem after this result.)

**Result 5.** *Conditional on $\{T > t\}$, the following holds.*

- **(Low intensity regime)** *For $\ell < \ell_{cr}$, the ball $B(0, \sqrt{2\beta}t)$ is emptied, the system stays inside this ball and branches at rate $\beta$.*
- **(High intensity regime)** *For $\ell > \ell_{cr}$,*
  - *$d = 1$: a $B(0, o(t))$ ball is emptied, while $Z$ suppresses branching until $t$ and stays inside this ball,*
  - *$d \geq 2$: There are two numbers (these are actually the unique minimizers of the variational problem for $I(\ell, \beta, d)$), $0 < \eta_* < 1$ and $c^* > 0$, such that $Z$ suppresses branching until $\eta^* t$, empties a ball of* radius $\sqrt{2\beta}(1-\eta^*)t$ around a point at distance $c^*t$ from the origin; during remaining time $(1-\eta^*)t$ the BBM branches at its normal $\beta$ rate.

Furthermore, for $d \geq 2$, $\ell \mapsto I(\ell, \beta, d)$ is continuous and strictly increasing, with $\lim_{\ell \to \infty} I(\ell, \beta, d) = \beta$ (see Fig. 1) and $\ell \mapsto \eta^*(\ell, \beta, d)$ and $\ell \mapsto c^*(\ell, \beta, d)$ are both discontinuous at $\ell_{cr}$ and continuous on $(\ell_{cr}, \infty)$, with

$$\lim_{\ell \to \infty} \frac{1 - \eta^*(\ell, \beta, d)}{c^*(\ell, \beta, d)} = 1, \qquad \lim_{\ell \to \infty} c^*(\ell, \beta, d) = 0. \qquad (20)$$

Finally, $c^* > \sqrt{2\beta}\,(1-\eta^*)$ for all $\ell > \ell_{cr}$.

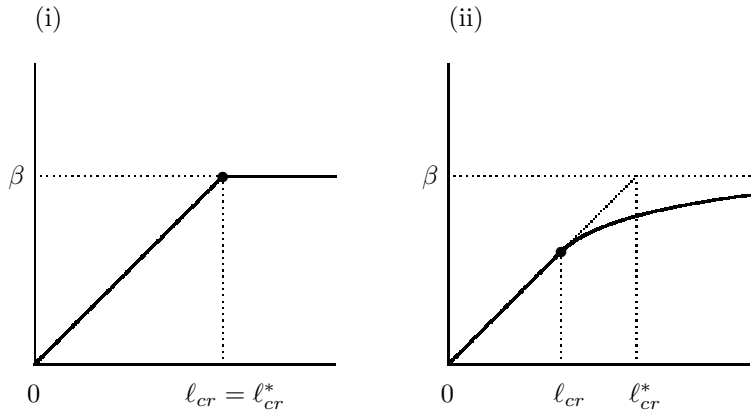

Fig. 1  $\ell \mapsto I(\ell, \beta, d)$ for: (i) $d = 1$; (ii) $d \geq 2$.

**Remark 10.** By the last statement, the clearing never contains the origin for $\ell > \ell_{cr}$. However, for $\ell < \ell_{cr}$, it does. So, the discontinuity mentioned in the theorem means that, surprisingly, when $\ell$ crosses its critical value, both $c^*$ (the distance of the clearing normalized by $t$) and $\sqrt{2\beta}(1-\eta^*)$ (the radius of the clearing normalized by $t$) have a jump. So in a sense, the clearing 'jumps away' from the origin and its size 'jumps down'. We have no intuitive explanation for this phenomenon.                                                                                           ◇

The precise statement regarding the optimal survival strategy is as follows.

**Theorem 10** (optimal survival strategy). *Fix $\beta, a$. For $r, b > 0$ and $t \geq 0$, define*

$$C(t; r, b) = \{\exists x_0 \in \mathbb{R}^d \colon \ |x_0| = b, B_{rt}(x_0 t) \cap K = \emptyset\}. \qquad (21)$$

(i) For $d = 1$, $\ell < \ell_{cr}$ or $d \geq 2$, any $\ell$, and $0 < \varepsilon < 1 - \eta^*$,

$$\lim_{t\to\infty} (\mathbb{E} \times P_{\delta_0})\Big(C\left(t; \sqrt{2\beta}\,(1 - \eta^* - \varepsilon), c^*\right) \mid T > t\Big) = 1,$$
$$\lim_{t\to\infty} (\mathbb{E} \times P_{\delta_0})\Big(|Z(t)| \geq \lfloor e^{\beta(1-\eta^*-\varepsilon)t} \rfloor \mid T > t\Big) = 1. \tag{22}$$

(ii) For $d \geq 1$, $\ell < \ell_{cr}$ and $\varepsilon > 0$,

$$\lim_{t\to\infty} (\mathbb{E} \times P_{\delta_0})\Big(B_{(1+\varepsilon)\sqrt{2\beta}\,t}(0) \cap K \neq \emptyset \mid T > t\Big) = 1,$$
$$\lim_{t\to\infty} (\mathbb{E} \times P_{\delta_0})\Big(R(t) \subseteq B_{(1+\varepsilon)\sqrt{2\beta}\,t}(0) \mid T > t\Big) = 1,$$
$$\lim_{t\to\infty} (\mathbb{E} \times P_{\delta_0})\Big(R(t) \not\subseteq B_{(1-\varepsilon)\sqrt{2\beta}\,t}(0) \mid T > t\Big) = 1. \tag{23}$$

(iii) For $d \geq 1$, $\ell > \ell_{cr}$ and $0 < \varepsilon < \eta^*$,

$$\lim_{t\to\infty} (\mathbb{E} \times P_{\delta_0})\Big(|Z((\eta^* - \varepsilon)t)| \leq \lfloor t^{d+\varepsilon} \rfloor \mid T > t\Big) = 1. \tag{24}$$

(iv) For $d = 1$, $\ell > \ell_{cr}$ and $\varepsilon > 0$,

$$\lim_{t\to\infty} (\mathbb{E} \times P_{\delta_0})\Big(B_{\varepsilon t}(0) \cap K \neq \emptyset \mid T > t\Big) = 1,$$
$$\lim_{t\to\infty} (\mathbb{E} \times P_{\delta_0})\Big(R(t) \subseteq B_{\varepsilon t}(0) \mid T > t\Big) = 1. \tag{25}$$

When $d = 1$ we get an answer for the one dimensional constant intensity case. In fact, the analysis shows that the optimal strategy is more extreme: in the low intensity regime $I(\ell, \beta, d)$ grows linearly with $\ell$ and reaches the value $\beta$, while in the high intensity regime it is always $\beta$ (see Fig. 1). This corresponds to the picture that in the low intensity regime there is neither suppressing of branching nor dislocation ($\eta^* = c^* = 0$), but an interval around the origin must be emptied, whereas in the high intensity regime, the branching is completely suppressed ($\eta^* = 1, c^* = 0$).

The following problems are open:

- For $d = 1$ and $\ell > \ell_{cr}$, what is the radius of the $o(t)$-ball that is emptied and how many particles are there inside this ball at time $t$?
- For $d \geq 2$ and $\ell > \ell_{cr}$, what is the shape of the "small tube" in which the system moves its particles away from the origin while suppressing the branching? How many particles are alive at time $\eta^* t$?
- What can be said about the optimal survival strategy at $\ell = \ell_{cr}$?
- Instead of letting the trap density decay to zero at infinity, another way to make survival easier is by giving the Brownian motion an inward drift while keeping the trap density field constant. Suppose that $d\nu/dx \equiv \ell$ and that the inward drift *radially increases* like $\sim \kappa |x|^{d-1}$, $|x| \to \infty$, $\kappa > 0$. Is there again a crossover in $\ell$ at some critical value $\ell_{cr} = \ell_{cr}(\kappa, \beta, d)$? What is the optimal survival strategy?

We close this section with discussing briefly another killing rule, which we call *'individual killing rule'*. We now only kill $Z$ when all its particles have been absorbed. So the death time becomes

$$\widetilde{T} := \inf\{t \geq 0 \mid |Z_t^K| = 0\},$$

where $Z^K$ is the BBM with killing at $\partial K$.

What can we say about the tail behavior of the survival probability then? Well, actually there is no tail behavior because there is no tail at all! And this is true for *any* locally finite trap intensity!

More precisely, the following holds.

**Theorem 11.** *Fix $d, \beta, a$. For any locally finite intensity measure $\nu$,*

$$\lim_{t\to\infty} (\mathbb{E} \times P_{\delta_0})(\tilde{T} > t) > 0. \tag{26}$$

Theorem 11 follows from the assertion that the system may survive by emptying a ball with a finite radius $R > R_0$, where $R_0$ is chosen such that the branching rate $\beta$ balances against the killing rate $\lambda(B(0, R_0))$, the principal Dirichlet eigenvalue of $-\Delta/2$ on the ball $B(0, R_0)$: $\lambda(B(0, R_0)) = \beta$. Indeed, it can be shown, that for any $R > R_0$ there is a strictly positive probability (denoted by $p_R$) that at all times at least one particle has not yet left $B(0, R_0)$. Consequently, the survival probability is bounded from below by $\sup_{R > R_0} \{p_R \exp[-\nu(B(0, R_0))]\}$.

The following two problems for individual killing are open:

– What is the limit in (26), say, when $d\nu/dx$ is spherically symmetric?
– If the Brownian motion is given an outward drift, then for what values of the drift does the survival probability decay to zero?

## 9. Spatial branching processes among Poissonian obstacles: 'mild' obstacles

### *9.1. Basic notions and questions*

Let us consider now a mechanism which is 'milder' than killing but still reduces the population size. The purpose is to study a spatial branching model with the property that the branching rate is decreased in a certain random region. More specifically we will use a natural model for the random environment: *mild Poissonian obstacles*.

We are going to work with slightly different notation relative to the previous sections. As before, let $\omega$ be a Poisson point process (PPP) on $\mathbb{R}^d$ with intensity $\nu > 0$ and let now **P** denote the corresponding law. Furthermore, let $a > 0$ and $0 < \beta_1 < \beta_2$ be fixed. We define the *branching Brownian motion (BBM) with a mild Poissonian obstacle*, or the '$(\nu, \beta_1, \beta_2, a)$-BBM' as follows. Let $B(z, r)$ denote the open ball centered at $z \in \mathbb{R}^d$ with radius $r$ and let $K$ denote the random set given by the $a$-neighborhood of $\omega$:

$$K = K_\omega := \bigcup_{x_i \in \mathrm{supp}(\omega)} \overline{B}(x_i, a).$$

Then $K$ is a *mild obstacle configuration* attached to $\omega$. This means that given $\omega$, we define $P^\omega$ as the law of the strictly dyadic (i.e. precisely two offspring) BBM on $\mathbb{R}^d$, $d \geq 1$ with spatially dependent branching rate

$$\beta(x, \omega) := \beta_1 \mathbb{1}_{K_\omega}(x) + \beta_2 \mathbb{1}_{K_\omega^c}(x).$$

An equivalent (informal) definition is that as long as a particle is in $K^c$, it obeys the branching rule with rate $\beta_2$, while in $K$ its reproduction is suppressed and it branches with the smaller rate $\beta_1$. (We assume that the process starts with a single particle at the origin.) The process under $P^\omega$ is called a *BBM with mild Poissonian obstacles* and denoted by $Z$. The total mass process will be denoted by $|Z|$. Further, $W$ will denote $d$-dimensional Brownian motion with probabilities $\{\mathbb{P}_x, \ x \in \mathbb{R}^d\}$.

**Remark 11** (self-duality). The discrete setting has the advantage that when $\mathbb{R}^d$ is replaced by $\mathbb{Z}^d$, the difference between the sets $K$ and $K^c$ is no longer relevant. Indeed the equivalent of a Poisson trap configuration is an i.i.d. trap configuration on the lattice, and then its complement is also i.i.d. (with a different parameter). So, in the discrete case 'Poissonian mild obstacles' give the same type of model as 'Poissonian mild catalysts' would. This nice duality is lost in the continuous setting as the 'Swiss cheese' $K^c$ is not the same type of geometric object as $K$. ◇

Consider now the following natural questions (both in the annealed and the quenched sense):

1. What can we say about the growth of the total population size?

2. What are the large deviations? (E.g., what is the probability of producing an atypically small population.)
3. What can we say about the *local* population growth?

As far as (1) is concerned, recall that the total population of an ordinary (free) BBM grows a.s. and in expectation as $e^{\beta_2 t}$. (Indeed, for ordinary BBM, the spatial component plays no role, and hence the total mass is just a $\beta_2$-rate pure birth process $X$. As is well known, the limit $N := \lim_{t\to\infty} e^{-\beta_2 t} X_t$ exists a.s and in mean, and $P(0 < N < \infty) = 1$.) In our model of BBM with the reproduction blocking mechanism, how much will the suppressed branching in $K$ slow the global reproduction down? Will it actually change the exponent $\beta_2$? (We will see that although the global reproduction does slow down, the slowdown is captured by a sub-exponential factor, being different for the quenched and the annealed case.)

Consider now (2). Here is an argument to give a further motivation. Let us assume for simplicity that $\beta_1 = 0$ and ask the simplest question: what is the probability that there is no branching at all up to time $t > 0$? In order to avoid branching the first particle has to 'resist' the branching rate $\beta_2$ inside $K^c$. Therefore this question is quite similar to the survival asymptotics for a single Brownian motion among 'soft obstacles' — but of course in order to prevent branching the particle seeks for large islands covered by $K$ rather then the usual 'clearings'. In other words, the particle now prefers to avoid the $K^c$ instead of $K$. Hence, (2) above is a possible generalization of this (modified) soft obstacle problem for a single particle, and the presence of branching creates new type of challenges.

Finally, for further mathematical models in population biology see e.g. [4].

### 9.2. Some results on the population growth

Let us consider the total population size first. Since the behavior of the total population size does not seem to be easy to handle, let us first try to handle its *expectation* and then trusting in the Law of Large Numbers (i.e. that the total population size behaves as its expectation), we will at least *know*, what we want to prove.

Concerning the expected global growth rate we have the following result, which, although not at all trivial, can very easily be derived from the Donsker-Varadhan result.

**Lemma 2** (Expected global growth rate). *On a set of full* **P***-measure,*

$$E^\omega |Z_t| = \exp\left[\beta_2 t - c(d,\nu)\frac{t}{(\log t)^{2/d}}(1+o(1))\right], \text{ as } t \to \infty \qquad (27)$$

*(quenched asymptotics), and*

$$(\mathbf{E} \otimes E^\omega)|Z_t| = \exp[\beta_2 t - \tilde{c}(d,\nu)t^{d/(d+2)}(1+o(1))], \text{ as } t \to \infty \qquad (28)$$

*(annealed asymptotics)*, where

$$c(d,\nu) := \lambda_d \left(\frac{d}{\nu \omega_d}\right)^{-2/d},$$

$$\tilde{c}(d,\nu) := (\nu \omega_d)^{2/(d+2)} \left(\frac{d+2}{2}\right) \left(\frac{2\lambda_d}{d}\right)^{d/(d+2)},$$

and $\omega_d$ is the volume of the $d$-dimensional unit ball, while $\lambda_d$ is the principal Dirichlet eigenvalue of $-\frac{1}{2}\Delta$ on it.

Notice that

1. $\beta_1$ does not appear in the formulas,
2. the higher the dimension is, the smaller the expected population size is.

**Remark 12.** Let us pretend for a moment that we are talking about an ordinary BBM. Then at time $t$ one has $e^{\beta_2 t}$ particles with probability tending to one as $t \to \infty$ (the population size divided by $e^{\beta_2 t}$ has a limit, to be precise). For $t$ fixed take a ball $B = B(0, R)$ (here $R = R(t)$) and let $K$ be so that $B \subset K^c$ (such a ball left empty by $K$ is called a *clearing*). Consider the expected number of particles that are confined to $B$ up to time $t$. These particles do not feel the blocking effect of $K$, while the other particles may have not been born due to it.

Optimize $R(t)$ with respect to the cost of having such a clearing and the probability of confining a single Brownian motion to it. This is precisely the same optimization as for the classical Wiener-sausage. Hence one gets the expectation in the theorem as a lower estimate.

One suspects that the main contribution in the expectation in (28) is coming from the expectation on the event of having a clearing with optimal radius $R(t)$. In other words, denoting by $p_t$ the probability that a single Brownian particle stays in the $R(t)$-ball up to time $t$, one argues heuristically that $p_t e^{\beta_2 t}$ particles will stay inside the clearing up to time $t$ 'for free' (i.e. with probability tending to one as $t \to \infty$).

The intuitive reasoning is as follows. If we had *independent* particles instead of BBM, then, by a 'Law of Large Numbers type argument' (using Chebysev inequality and the fact that $\lim_{t\to\infty} p_t e^{\beta_2} = \infty$), roughly $p_t e^{\beta_2 t}$ particles out of the total $e^{\beta_2 t}$ would stay in the $R(t)$-ball up to time $t$ with probability tending to 1 as $t \uparrow \infty$. One suspects then that the lower estimate remains valid for the branching system too, because the particles are 'not too much correlated'. This kind of argument (in the quenched case though) can be made precise by estimating certain covariances.

To understand *the difference between the annealed and the quenched case* note that in the latter case large clearings (far away) are automatically (that is, **P**-a.s.) present. Hence, similarly to the single Brownian particle problem, the difference between the two asymptotics is due to the fact that even though there is an appropriate clearing far away **P**-a.s., there is one around the origin with a small (but not too small) probability. Still, the two cases will have a similar element when one drops the expectations and investigate the process itself, and show that inside such a clearing a large population is going to flourish. This also leads to the intuitive explanation for the *decrease of the population size as the dimension increases.* The radial drift in dimension $d$ is $(d-1)/2$. The more transient the motion is (that is, the larger $d$ is), the harder it is for the particles to stay in the appropriate clearings. ◇

Let us prove now the expectation result. First we need a lemma.

**Lemma 3** (Expectation given by Brownian functional)**.** *Fix $\omega$. Then*

$$E^\omega |Z_t| = \mathbb{E} \exp\left[\int_0^t \beta(W_s)\right] \mathrm{d}s. \tag{29}$$

**Proof.** It is well known ('first moment formula' of spatial branching processes), that $E_x^\omega |Z_t| = (T_t 1)(x)$, where $u(x,t) := (T_t 1)(x)$ is the minimal solution of the parabolic problem:

$$\begin{aligned}
\frac{\partial u}{\partial t} &= \left(\frac{1}{2}\Delta + \beta\right) u \text{ on } \mathbb{R}^d \times (0, \infty), \\
u(\cdot, 0) &= 1, \\
u &\geq 0.
\end{aligned} \tag{30}$$

(Here $\{T_t\}_{t \geq 0}$ denotes of course the semigroup corresponding to the generator $\frac{1}{2}\Delta + \beta$ on $\mathbb{R}^d$.) This is equivalent (by the Feynman-Kac formula) to (29). □

**Proof of Lemma 2.**: Since $\beta := \beta_1 \mathbb{1}_K + \beta_2 \mathbb{1}_{K^c} = \beta_2 - (\beta_2 - \beta_1)\mathbb{1}_K$, we can rewrite the equation (29) as

$$E^\omega |Z_t| = e^{\beta_2 t} \mathbb{E} \exp\left[-\int_0^t (\beta_2 - \beta_1)\mathbb{1}_K(W_s)\,\mathrm{d}s\right].$$

The expectation on the righthand side is precisely the survival probability among soft obstacles of 'height' $\beta_2 - \beta_1$. (As we have already mentioned, that one does not sum the shape functions on the overlapping balls, does not make any difference with regard to the asymptotics) . The statements thus follow from the Donsker-Varadhan asymptotics for soft obstacles . □

We now investigate the behavior of the (quenched) global growth rate.

As already mentioned in these notes, it is a notoriously hard problem to prove the Law of Large Numbers in full generality for spatial branching systems, and the not purely exponential case is particularly challenging.

To elucidate this point, let us consider the $(L, \beta, D)$-branching diffusion with some Euclidean domain $D \subset \mathbb{R}^d$. Let $0 \not\equiv f$ be a nonnegative compactly supported bounded measurable function on $D$. If $\lambda_c$, the generalized principal eigenvalue of $L + \beta$ on $D$ is positive and $S = \{S_t\}_{t \geq 0}$ denotes the semigroup corresponding to $L+\beta$ on $D$, then $(S_t f)(\cdot)$ grows (pointwise) as $e^{\lambda_c t}$ on an exponential scale. However, in general, the scaling is not precisely exponential due to the presence of a subexponential term.

In the present case it turns out that the resulting operator $\frac{1}{2}\Delta + \beta$ does not scale precisely exponentially **P**-a.s. Replacing $f$ by the function $g \equiv 1$, it is still true that the growth is not precisely exponential – this is readily seen in Claim 2 and its proof.

Since the process in expectation is related to the semigroup $S$, therefore purely exponential scaling indicates the same type of scaling for the expectation of the process (locally or globally). Although, if one is interested in the scaling of the process itself (not just the expectation), the case when there is an additional subexponential factor is much harder, in our model the *randomization* of the branching rate $\beta$ helps, as $\beta$ has some 'nice' properties for almost every environment $\omega$, i.e. the 'irregular' branching rates 'sit in the **P**-zero set'.

Define the *average growth rate* by

$$r_t = r_t(\omega) := \frac{\log |Z_t(\omega)|}{t}.$$

Replace now $|Z_t(\omega)|$ by its expectation $\overline{Z}_t := E^\omega |Z_t(\omega)|$ and define

$$\widehat{r}_t = \widehat{r}_t(\omega) := \frac{\log \overline{Z}_t}{t}.$$

Recall from Lemma 2, that on a set of full **P**-measure,

$$\lim_{t \to \infty} (\log t)^{2/d} (\widehat{r}_t - \beta_2) = -c(d, \nu). \tag{31}$$

It turns out that an analogous statement holds for $r_t$ itself.

**Theorem 12** (LLN). *On a set of full **P**-measure,*

$$\lim_{t \to \infty} (\log t)^{2/d} (r_t - \beta_2) = -c(d, \nu) \qquad \text{in } P^\omega - \text{probability}. \tag{32}$$

One interprets Theorem 12 as a kind of quenched law of large numbers. Loosely speaking,

$$r_t \approx \beta_2 - c(d, \nu)(\log t)^{-2/d} \approx \widehat{r}_t, \qquad t \to \infty,$$

on a set of full **P**-measure.

While the lower estimate is hard, the upper estimate follows easily from the expectation result. We only show here the latter one, but we will explain the 'philosophy' for the lower estimate too.

**Upper estimate.** Let $\epsilon > 0$. Using the Markov inequality along with the expectation formula (27), we have that on a set of full **P**-measure:

$$P^\omega \left[ (\log t)^{2/d}(r_t - \beta_2) + c(d,\nu) > \epsilon \right]$$
$$= P^\omega \left\{ |Z_t| > \exp\left[ t\left( \beta_2 - c(d,\nu)(\log t)^{-2/d} + \epsilon(\log t)^{-2/d} \right) \right] \right\}$$
$$\leq E^\omega |Z_t| \cdot \left( \exp\left[ t\left( \beta_2 - c(d,\nu)(\log t)^{-2/d} + \epsilon(\log t)^{-2/d} \right) \right] \right)^{-1}$$
$$= \exp\left[ -\epsilon t(\log t)^{-2/d} + o\left( t(\log t)^{-2/d} \right) \right] \to 0, \qquad \text{as } t \to \infty. \qquad \square$$

**Lower estimate (sketch).** Let us just outline the *strategy* of the proof. A key step is introducing three different time scales, $\ell(t)$, $m(t)$ and $t$ where $\ell(t) = o(m(t))$ and $m(t) = o(t)$ as $t \to \infty$. For the first, shortest time interval, one uses that there are 'many' particles produced and they are not moving 'too far away', for the second (of length $m(t) - \ell(t)$) one uses that one particle moves into a clearing of a certain size at a certain distance, and in the third one (of length $t - m(t)$) one uses that there is a branching tree emanating from that particle so that a certain proportion of particles of that tree stay in the clearing with probability tending to one.

It turns out that for example the following choices of $\ell$ and $m$ are appropriate: let $\ell(t)$ and $m(t)$ be arbitrarily defined for $t \in [0, e]$, and

$$\ell(t) := t^{1-1/(\log \log t)}, \ m(t) := t^{1-1/(2 \log \log t)}, \text{ for } t \geq t_0 > e.$$

Also, regarding the Poissonian environment one needs the following fact: let

$$R_0 = R_0(d, \nu) := \sqrt{\frac{\lambda_d}{c(d, \nu)}} = \left(\frac{d}{\nu \omega_d}\right)^{1/d},$$

(where $\lambda_d$ is the principal Dirichlet eigenvalue of $-\frac{1}{2}\Delta$ on the $d$-dimensional unit ball, and $\omega_d$ is the volume of that ball) and let

$$\rho(\ell) := R_0 (\log \ell)^{1/d} - (\log \log \ell)^2, \ t \geq 0.$$

Then,

$\mathbf{P}(\exists \ \ell_0(\omega) > 0 \ \text{ such that } \ \forall \ell > \ell_0(\omega) \ \exists \text{ clearing } B(x_0, \rho(\ell)) \text{ with } |x_0| \leq \ell) = 1.$

The full proof is beyond the scope of these notes.

### 9.3. Some results on the spatial spread

One may wonder *how much the speed (spatial spread) of free BBM reduces* due to the presence of the mild obstacle configuration. Note that we are not talking about the bulk of the population (or the 'shape') but rather about *individual* particles traveling to very large distances from the origin.

As is well known, ordinary 'free' branching Brownian motion with constant branching rate $\beta_2 > 0$ has radial speed $\sqrt{2\beta_2}$. Let $N_t$ denote the population size at $t \geq 0$ and let $\xi_k$ ($1 \leq k \leq N_t$) denote the position of the $k$th particle (with arbitrary labeling) in the population. Furthermore, let $m(t)$ denote a number for which $u(t, m(t)) = \frac{1}{2}$, where

$$u(t, x) := P\left[\max_{1 \leq k \leq N_t} \|\xi_k(t)\| \leq x\right].$$

In a classic paper, Bramson considered the one dimensional case and proved that as $t \to \infty$,

$$m(t) = t\sqrt{2\beta_2} - \frac{3}{2\sqrt{2\beta_2}} \log t + \mathcal{O}(1). \tag{33}$$

Since the one-dimensional projection of a $d$-dimensional branching Brownian motion is a one-dimensional branching Brownian motion, we can utilize Bramson's result for the higher dimensional cases too. Namely, it is clear, that in high dimension the spread is *at least* as quick as in (33). The asymptotics (33) is derived for the case $\beta_2 = 1$; the general result can be obtained similarly.

Studying the function $u$ has significance in analysis too as $u$ solves

$$\frac{\partial u}{\partial t} = \frac{1}{2} u_{xx} + \beta_2 (u^2 - u), \tag{34}$$

with initial condition

$$u(0, \cdot) = \mathbb{1}_{[0,\infty)}(\cdot). \tag{35}$$

We will see below that the branching Brownian motion with mild obstacles *spreads less quickly* than ordinary branching Brownian motion by giving an upper estimate on its speed.

**Remark 13.** A related result was obtained earlier by Lee-Torcaso [10], but, unlike in (33), only up to the linear term and moreover, for random walks instead of Brownian motions. Their approach was to consider the problem as the description of wave-front propagation for a random KPP equation.    ⋄

Before turning to the upper estimate, we discuss the lower estimate. We are going to show that, if in our model Brownian motion is replaced by Brownian motion with constant drift $\gamma$ in a given direction, then any fixed nonempty ball is recharged infinitely often with positive probability, as long as the drift satisfies $|\gamma| < \sqrt{2\beta_2}$.

For simplicity, assume that $d = 1$ (the general case is similar). Fix the environment $\omega$. Recall Doob's $h$-transform of second order elliptic operators:

$$L^h(\cdot) := \frac{1}{h} L(h \cdot).$$

Applying an $h$-transform with $h(x) := \exp(-\gamma x)$, a straightforward computation shows that the operator

$$L := \frac{1}{2} \frac{d^2}{dx^2} + \gamma \frac{d}{dx} + \beta_2$$

transforms into

$$L^h = \frac{1}{2} \frac{d^2}{dx^2} - \frac{\gamma^2}{2} + \beta_2.$$

Then one can show that the generalized principal eigenvalue for this latter operator is $-\frac{\gamma^2}{2} + \beta_2$ for almost every environment. Since the generalized principal eigenvalue is invariant under $h$-transforms, it follows that $-\frac{\gamma^2}{2} + \beta_2 > 0$ is the generalized principal eigenvalue of $L$. By subsection 3.2 then, any fixed nonempty interval is recharged infinitely often with positive probability.

Turning back to our original setting, the application of the '*spine*'-technology seems also promising. As we have seen in section 3., the spine method uses a Girsanov-type change of measure on the branching diffusion. The question of uniform integrability for the martingale density can be checked and under the new measure, the process can be represented 'directly' by a 'spine' type construction. The 'spine' represents a special '*deviant*' particle performing a motion which has important properties from the point of view of the problem regarding the original branching process. The question regarding 'going back to the original measure' is then related to the UI property of the martingale.

In our case perhaps it is not very difficult to show the existence of a 'spine' particle (under a martingale-change of measure) that has drift $\gamma$ as long as $|\gamma| < \sqrt{2\beta_2}$.

We close this section with an upper estimate (without proof) in which the order of the correction term is *larger* than the $\mathcal{O}(\log t)$ term appearing in Bramson's result, namely it is $\mathcal{O}\left(\frac{t}{(\log t)^{2/d}}\right)$. (All orders are meant for $t \to \infty$.) The following theorem says that, whatever $\beta_1 \in (0, \beta_2)$ is, loosely speaking, at time $t$ the spread is not more than

$$t\sqrt{2\beta_2} - c(d, \nu)\sqrt{\frac{\beta_2}{2}} \cdot \frac{t}{(\log t)^{2/d}}.$$

Recall the notation: $c(d, \nu) := \lambda_d \left(\frac{\nu\omega_d}{d}\right)^{2/d}$, $\omega_d$ is the volume of the unit ball in $\mathbb{R}^d$, and $\lambda_d$ is the principal Dirichlet eigenvalue of $-\frac{1}{2}\Delta$ on it.

**Theorem 13.** *Define the functions*

$$f(t) := c(d, \nu)\frac{t}{(\log t)^{2/d}} \quad \text{and} \quad n(t) := t\sqrt{2\beta_2} \cdot \sqrt{1 - \frac{f(t)}{\beta_2 t}}.$$

*Then*

$$n(t) = t\sqrt{2\beta_2} - c(d, \nu)\sqrt{\frac{\beta_2}{2}} \cdot \frac{t}{(\log t)^{2/d}} + \mathcal{O}\left(\frac{t}{(\log t)^{4/d}}\right); \text{ and if}$$

$$\begin{aligned} A_t &:= \quad \{\text{no particle has left the n(t) − ball up to t}\} \\ &= \quad \left\{\bigcup_{0 \le s \le t} \text{supp}(Z_s) \subseteq B(0, n(t))\right\}, \end{aligned}$$

*then on a set of full* **P***-measure,*

$$\liminf_{t \to \infty} P^\omega(A_t) > 0.$$

### 9.4. Branching L-diffusion with mild obstacles

We now generalize the setting for the case when the underlying motion is a diffusion. Let **P** be as before but replace the Brownian motion by an $L$-diffusion $Y$ on $\mathbb{R}^d$, where $L$ is a second order elliptic operator as in Subsection 2.1. The branching $L$-diffusion with the Poissonian obstacles can be defined analogously to the case of Brownian motion.

The following result demonstrates that the local behavior of the process exhibits a *dichotomy*. The crossover is given in terms of the local branching rate $\beta_2$ and $\lambda_c(L)$: as already explained in Subsection 2.1, local extinction occurs when the branching rate inside the 'free region' $K^c$ is not sufficiently large to compensate the transience of the underlying $L$-diffusion; if it is strong enough, then local mass grows exponentially.

There is, however, an interesting new feature of the result: *neither $\beta_1$ nor the intensity $\nu$ of the obstacles play role* .

**Claim 1** (Quenched exponential growth/local extinction). *Given the environment $\omega$, denote by $P^\omega$ the law of the branching $L$-diffusion.*

(i) *Let $\beta_2 > -\lambda_c(L)$ and let $\nu > 0$ and $\beta_1 \in (0, \beta_2)$ be arbitrary. Then the following holds on a set of full $\mathbf{P}$-measure: For any $\epsilon > 0$ and any bounded open set $\emptyset \neq B \subset \mathbb{R}^d$,*

$$P^\omega \left( \limsup_{t \uparrow \infty} e^{(-\beta_2 - \lambda_c(L) + \epsilon)t} Z_t(B) = \infty \right) > 0.$$

*and*

$$P^\omega \left( \limsup_{t \uparrow \infty} e^{(-\beta_2 - \lambda_c(L))t} Z_t(B) < \infty \right) = 1.$$

(ii) *Let $\beta_2 \leq -\lambda_c(L)$ and let $\nu > 0$ and $\beta_1 \in (0, \beta_2)$ be arbitrary. Then the following holds on a set of full $\mathbf{P}$-measure: For any bounded open set $B \subset \mathbb{R}^d$ there exists a $P^\omega$-a.s. finite random time $t_0 = t_0(B)$ such that $X_t(B) = 0$ for all $t \geq t_0$, (local extinction).*

*Proof.* In order to be able to use Theorem 3, we compare the rate $\beta$ with another, smooth (i.e. $C^\gamma$) function $V$. Recalling that $K = K_\omega := \bigcup_{x_i \in \text{supp}(\omega)} \overline{B}(x_i, a)$, let us enlarge the obstacles:

$$K^* = K^*_\omega := \bigcup_{x_i \in \text{supp}(\omega)} \overline{B}(x_i, 2a).$$

Then $(K^*)^c \subset K^c$. Recall that $\beta(x) := \beta_1 \mathbb{1}_K(x) + \beta_2 \mathbb{1}_{K^c}(x) \leq \beta_2$ and let $V \in C^\gamma$ ($\gamma \in (0, 1]$) with[16]

$$\beta_2 \mathbb{1}_{(K^*)^c} \leq V \leq \beta. \tag{36}$$

In fact it is easy to see the existence of such functions which are even $C^\infty$ by writing $\beta = \beta_2 - (\beta_2 - \beta_1) \mathbb{1}_K$ and considering the function $V := \beta_2 - (\beta_2 - \beta_1)f$, where $f \geq \mathbb{1}_{K^*}$ and $f$ is a $C^\infty$-function obtained as follows. Let $f$ be a sum of compactly supported $C^\infty$-functions $f_n$, $n \geq 1$, with disjoint support, where

---
[16]The existence of a *continuous* function satisfying (36) would of course immediately follow from Uryson's Lemma.

supp($f_n$) is in the $\epsilon_n$-neighborhood of the $n$th connected component of $K^*$, with appropriately small $0 < \epsilon_n$'s. Consider the operator $L + V$ on $\mathbb{R}^d$ and let $\lambda_c = \lambda_c(\omega)$ denote its generalized principal eigenvalue. Since $V \in C^\gamma$, we are in the setting of Theorem 3, and furthermore, since $V \leq \beta_2$, $\lambda_c \leq \lambda_c(L) + \beta_2$ for every $\omega$.

On the other hand, one gets a lower estimate on $\lambda_c$ as follows. Fix $R > 0$. Since $\beta_2 \mathbb{1}_{(K^*)^c} \leq V$, by the homogeneity of the Poisson point process, for almost every environment the set $\{x \in \mathbb{R}^d \mid V(x) = \beta_2\}$ contains a clearing of radius $R$. Hence, by comparison, $\lambda_c \geq \lambda^{(R)}$, where $\lambda^{(R)}$ is the principal Dirichlet eigenvalue of $L + \beta_2$ on a ball of radius $R$. Since $R$ can be chosen arbitrarily large and since $\lim_{R \uparrow \infty} \lambda^{(R)} = \lambda_c(L) + \beta_2$, we conclude that $\lambda_c \geq \lambda_c(L) + \beta_2$ for almost every environment.

From the lower and upper estimates, we obtain that

$$\lambda_c = \lambda_c(L) + \beta_2 \text{ for a.e. } \omega. \tag{37}$$

Consider now the branching processes with the same motion component $L$ but with rate $V$, respectively constant rate $\beta_2$. The statements (i) and (ii) of Claim 1 are true for these two processes by (37) and Theorem 3. As far as the original process (with rate $\beta$) is concerned, (i) and (ii) of Claim 1 now follow by $\omega$-wise comparison. □

## 10. Generalizations and open problems for mild obstacles

In this section we suggest some further problems and directions for research.

### *10.1. More general branching*

It should also be investigated, what happens when the dyadic branching law is replaced by a general one (but the random branching rate is as before). In a more sophisticated population model, particles can also die — then the obstacles do not necessarily reduce the population size as they sometimes prevent death.

#### *10.1.1. Supercritical branching*

When the offspring distribution is supercritical, the method of our paper seems to work, although when the offspring number can also be zero, one has to *condition on survival* for getting the asymptotic behavior.

#### *10.1.2. (Sub)critical branching*

Critical branching requires an approach very different from the supercritical one, since taking expectations now does not provide a clue: $E^\omega |Z_t(\omega)| = 1$, $\forall t > 0$, $\forall \omega \in \Omega$.

Having the obstacles, the first question is whether it is still true that

$$P^\omega(\text{extinction}) = 1 \ \forall \omega \in \Omega.$$

The answer is yes. To see this, note that since $|Z|$ is still a martingale, it has a nonnegative a.s. limit. This limit must be zero; otherwise $|Z|$ would stabilize at a positive integer. This, however is impossible because following one Brownian particle it is obvious that this particle experiences branching events for arbitrarily large times.

Setting $\beta_1 = 0$, the previous argument still goes through for almost all environments, because for almost all environments, the generic Brownian particle will eventually branch, and so in fact it will branch infinitely often. Let $\tau$ denote the almost surely finite extinction time for this case. One of the basic questions is the decay rate for $P^\omega(\tau > t)$. Will the tail be significantly heavier than $\mathcal{O}(1/t)$? [Of course $1/t$ would be the rate without obstacles.]

The subcritical case can be treated in a similar fashion. In particular, the total mass is a supermartingale and $P^\omega(\text{extinction}) = 1 \ \forall \omega \in \Omega$.

### *10.2. Superprocesses with mild obstacles*

An alternative view on the BBM with mild obstacles is as follows. Arguably, the model can be viewed as a *catalytic* BBM as well — the catalytic set is then $K^c$ (in the sense that branching is 'intensified' there). Catalytic spatial branching (mostly for superprocesses though) has been the subject of vigorous research in the last twenty years initiated by Dawson, Fleischmann and others. In those models the individual branching rates of particles moving in space depend on the amount of contact between the particle ('reactant') and a certain random medium called the catalyst. The random medium is usually assumed to be a 'thin' random set (that could even be just one point) or another superprocess. Sometimes 'mutually' or even 'cyclically' catalytic branching is considered. (See the next section for more on catalytic branching.)

Although there is a large amount of ongoing research on catalytic superprocesses, in those models, one usually cannot derive sharp quantitative results.

In light of this connection with catalytic superprocesses, a further goal is to generalize our setting by defining *superprocesses with mild obstacles* analogously to the BBM with mild obstacles.

Consider the $(L, \beta, \alpha, \mathbb{R}^d)$–superdiffusion, $X$. The definition of the superprocess with mild obstacles is straightforward: the parameter $\alpha$ on the (random) set $K$ is smaller than elsewhere.

Similarly, one can consider the case when instead of $\alpha$, the 'mass creation term' $\beta$ is random, for example with $\beta$ defined in the same way (or with a mollified version) as for the discreet branching particle system. Denote now by $P^\omega$ the law of this latter superprocess for given environment. We suspect that the superprocess with mild obstacles behaves similarly to the discrete branching process with mild obstacles when $\lambda_c(L + \beta) > 0$ and $P^\omega(\cdot)$ is replaced by $P^\omega(\cdot \mid X \text{ survives})$. The upper estimate can be carried out in a manner similar to the discrete particle system, as the expectation formula is still in force for superprocesses.

### 10.3. Unique dominant branch

Having discussed the growth rate of the population, the next step can be as follows. Once one knows the global population size $|Z_t|$, the model can be rescaled (normalized) by $|Z_t|$, giving a population of fixed weight. In other words, one considers the discrete probability measure valued process

$$\widetilde{Z}_t(\cdot) := \frac{Z_t(\cdot)}{|Z_t|}.$$

Then the question of the *shape* of the population for $Z$ for large times is given by the limiting behavior of the random probability measures $\tilde{Z}_t$, $t \geq 0$. (Of course, not only the particle mass has to be scaled, but also the spatial scales are interesting — see last paragraph.)

Can one for example locate a *unique dominant branch* (or, even if it is not unique, there are perhaps just a 'few' of them) for almost every environment, so that the total weight of its complement tends to (as $t \to \infty$) zero?

The motivation for this question comes from the proof of the lower estimate for Theorem 12. It seems conceivable that for large times the 'bulk' of the population will live in a clearing(s) within distance $\ell(t)$ and with radius

$$\rho(t) := R_0[\log \ell(t)]^{1/d} - [\log \log \ell(t)]^2.$$

### *10.4. Strong Law*

Last but not least, one cannot be fully satisfied with the convergence in probability in Theorem 3. To prove a.s. convergence (which, if true, is a kind of SLLN) seems to be rather challenging. In particular, the easy upper estimate does not work any longer.

## 11. Some catalytic branching problems

### 11.1. A discrete model

Let us start with a discrete catalytic model studied by H. Kesten and V. Sidoravicius. Here the branching particle system on $\mathbb{Z}^d$ is so that its branching is catalyzed by another autonomous particle system on $\mathbb{Z}^d$. There are two types of particles, the $A$-particles ('catalyst') and the $B$-particles ('reactant'). They move, branch and interact in the following way. Let $N_A(x,s)$ and $N_B(x,s)$ denote the number of $A$- [resp. $B$-]particles at $x \in \mathbb{Z}^d$ and at time $s \in [0,\infty)$. (All $N_A(x,0)$ and $N_B(x,0)$, $x \in \mathbb{Z}^d$ are independent Poisson variables with mean $\mu_A$ ($\mu_B$).)

Every $A$-particle ($B$-particle) performs independently a continuous-time random walk with jump rate $D_A$ ($D_B$). In addition a $B$-particle dies at rate $r$, and, when present at $x$ at time $s$, it splits into two in the next $ds$ time with probability $\beta N_A(x,s)ds + o(ds)$. Conditionally on the system of the $A$-particles, the jumps, deaths and splitting of the $B$-particles are independent.

Kesten and Sidoravicius showed for instance that for large $\beta$ there exists a critical $r_c$ such that it separates local extinction regime ($r \geq r_c$) from local survival regime ($r < r_c$). They also explained why the B-particles may survive if $\beta$ is large enough compared to the other parameters, even if one starts with only one A-particle and one B-particle in the entire system.

### 11.2. Catalytic superprocesses

The definition of catalytic *superprocesses* requires more care. The simplest case was introduced by D. Dawson and K. Fleischmann who investigated a one dimensional critical super-Brownian motion with a *single point catalyst*. This means formally that we set $\beta = 0$ and $\alpha = \delta_0$ (the Dirac-delta at zero). Intuitively this means that branching occurs only at a single point but there it occurs at an 'infinitely high' rate[17]. The precise definition is given through the log-Laplace equation, in which case (4) is replaced by an appropriate integral equation.

One of the interesting properties of this process is that if $P$ denotes the corresponding probability, then

$$P_\mu(|X_t| > 0, \ \forall t > 0, \ \text{but} \ |X_t| \to 0 \ \text{as} \ t \to \infty) = 1$$

for all $0 \neq \mu \in \mathcal{M}_f$, and $X_t(B) \to 0$ in $P_\mu$-probability for all $B \subset\subset \mathbb{R}$, even if $\mu$ is the Lebesgue measure. (The definition of a superprocess starting with a $\sigma$-finite measure is straightforward: it is the independent sum of superprocesses starting with finite measures.)

Another interesting case is, when, instead of the singular coefficient $\delta_0$ in the branching mechanism, one imagines that the catalytic set, where branching occurs, is given by the support of *another superprocess* (in low dimensions). Note that this is already the case of a *random catalyst*, that is, we are talking about a superprocess in a random medium (but, unlike in the cases we discussed, the medium *evolves in time*.) Since branching is critical, it is detrimental for the reactant to spend too much time in catalytic areas.

---

[17]Recall that since $\beta \equiv 0$, one can now identify $\alpha$ with the 'clock'.

For $d \leq 3$, Dawson and Fleischmann investigated the qualitative properties of this superprocess. They showed for example, that in dimension one, if both the catalyst and the reactant starts from Lebesgue measure then the reactant process survives locally.[18] This is dramatically different behavior from that of a classical one dimensional super-Brownian motion, since the latter would suffer local extinction even starting from $\mu =$ Lebesgue.

In two dimension, again starting both processes from Lebesgue measure, Dawson and Fleischmann proved that although the catalyst dies out locally, this is true only in probability. In particular, there is not a finite time after which a given ball becomes and remains empty. Consequently, at late times $T$, huge clusters of the catalyst (their height turns out to be of order $\log T$) come back to any given finite window in $\mathbb{R}^2$ and (since critical binary branching with infinite rate degenerates to pure killing) could kill all the (recurrent) reactant particles there. Thus it is somewhat surprising that in fact, the reactant *survives* locally.

Now let us go one step further, and imagine that we have two one dimensional, critical super-Brownian motions and they catalyze *each other*. That is, the first process only branches in the presence of the second process, and vice versa ('*mutually catalytic branching*'). Then, of course it becomes a nontrivial business to rigorously prove that such a model is well defined; in fact the log-Laplace equation approach breaks down. Such a model was first studied by D. Dawson and E. Perkins [7]. The definition uses a system of stochastic partial differential equations (SPDE's) for the densities of the two superprocesses, which by a result of L. Mytnik [19], uniquely determine the model.

Further models, including the *cyclically catalytic* setting, are beyond the scope of these lecture notes. We suffice here by merely defining the latter model. Cyclically catalytic super-Brownian motion is a spatial branching process with $K$ types of populations involved ($K \geq 2$). Each process $X^k$ is (catalytic) super-Brownian motion in one dimension but with a branching rate that is governed by the density of the succeeding type $X^{k+1}$. Hence the $K$ single processes interact in a cyclic fashion:

$$0 \to K-1 \to K-2 \to \cdots \to 1 \to 0.$$

The literature of catalytic spatial branching processes and superprocesses is huge and we have just glimpsed at it. The interested reader is referred to the survey paper [16].

---

[18]In fact, they showed that for almost all realizations of the catalytic medium (starting with Lebesgue measure), the reactant at $t$ converges (stochastically) to the starting Lebesgue measure.

**12. Concluding remarks**

Spatial branching processes is a huge area. So is spatial random media. Although we tried to cull some interesting problems, there is no doubt that many papers that should have been mentioned were unfortunately skipped.

As far as superprocesses are concerned, it is impossible even to sketch the research that has been carried out by Aldous, Dawson, Delmas, Donelly, Dynkin, Evans, Fleischmann, Klenke, Kurtz, Le Gall, Le Jan, Mörters, Mytnik, Mueller, Perkins, Pinsky and others. (And for discrete branching processes, there are those work of Athreya, Biggins, S. Harris, T. Harris, Jagers, Kyprianou, Ney, Pemantle, Wakolbinger, and the whole French school, to name just a few...)

Similarly, we did not even attempt to describe all the contributions made to the field of spatial random media by scholars like Bolthausen, Comets, Cranston, Donsker, Gärtner, den Hollander, van den Berg, Molchanov, Sznitman, Varadhan, Zeitouni, Zerner and others.

My hope is that those whose work was not mentioned in these notes will (forgive the author and) send me their critical comments.